\documentclass[12pt,a4paper]{article}
\usepackage{amsmath,amsfonts,amssymb,amscd}
\usepackage{amsmath}
\usepackage{tikz-cd}
\usetikzlibrary{matrix, calc, arrows}

\usepackage[body={6.6in, 9.2in},left=1in, top=1.2in]{geometry}
\def\rk{\operatorname{rk}}
\def\H{\operatorname{H}}
\def\gr{\operatorname{gr}}
\def\Ker{\operatorname{Ker}}

\def\GL{\operatorname{GL}}

\def\pr{\operatorname{pr}}

\def\Mat{\operatorname{Mat}}

\newtheorem{theorem}{Theorem}
\newtheorem{lemma}[theorem]{Lemma}
\newtheorem{definition}[theorem]{Definition}
\newtheorem{corollary}[theorem]{Corollary}
\newtheorem{proposition}[theorem]{Proposition}

\usepackage[all]{xy}
\usepackage[active]{srcltx}
\usepackage[parfill]{parskip}

\newcounter{th}

\newcounter{le}

\newcounter{lem}

\newcounter{de}

\newcounter{ex}

\begin{document}

\begin{center}
 {\bf   \LARGE{Rigidity of flag supermanifolds}}
\end{center}
\bigskip
\bigskip

\begin{center}
     {\large Elizaveta Vishnyakova}
\end{center}

\bigskip

\bigskip

\begin{abstract}
We prove that under certain assumptions a supermanifold of flags is rigid, this is its complex structure does not admit any non-trivial small deformation. Moreover under the same assumptions we show that a supermanifold of flags is unique non-split supermanifold with given retract.

\end{abstract}

\bigskip

\section{ Introduction}

It is a classical result that any flag manifold is rigid, see \cite{Bott}. In other words its complex structure does not possess any non-trivial small deformation. In general this statement is false for a flag supermanifold, see \cite{BunOni,Vaintrob}. For instance the projective superspace $\mathbb {CP}^{1|m}$, where $m \geq 4$, see \cite{Vaintrob}, and the super-grassmannian $\mathbf
{Gr}_{2|2,1|1}$, see \cite{BunOni}, are not rigid.  In \cite{onirig} it was proved that the super-grassmannian  $\mathbf
{Gr}_{m|n,k|l}$ is rigid if $m,n,k,l$ satisfy the following conditions
\begin{equation}
\begin{split}
&0<k<m, \quad 0<l<n, \\
&(k,l)\ne (1,n-1), \,\, (m-1,1), (1,n-2),\,\, (m-2,1), \,\, (2,n-1), \,\,(m-1,2).
\end{split}
\end{equation}
The idea of the paper \cite{onirig} is to compute $1$-cohomology with values in the tangent sheaf showing their triviality. This implies the rigidity of the super-grassmannian in this case. 

In this paper we compute  $1$-cohomology with values in the tangent sheaf $\mathcal T$ of a flag supermanifold showing their triviality. Therefore we prove that under certain conditions any flag supermanifold is rigid. We use the results of \cite{onirig} and the fact that a supermanifold of flags of length $n$ is a superbundle with base space a super-grassmannian and fiber a flag supermanifold of length $n-1$.

\textbf{Acknowledgements:} 
E.~V. was partially  supported by FAPEMIG, grant APQ-01999-18, and by Tomsk State University, Competitiveness Improvement Program.

\section{Main definitions}

\subsection{Supermanifolds}

 We use the word ``supermanifold'' in the sense of Berezin and
 Leites, see \cite{BL,ley,Man} for detail.  Throughout the paper, we will be interested in the complex-analytic version of the theory.  More precisely, a {\it complex-analytic superdomain of dimension $n|m$} is a $\mathbb{Z}_2$-graded
 ringed space
 $ \mathcal{U} = (\mathcal{U}_0, \mathcal{F}_{\mathcal{U}_0} \otimes \bigwedge(m) ),$
 where $\mathcal{F}_{\mathcal{U}_0}$ is the sheaf of holomorphic functions on an open set $\mathcal{U}_0\subset \mathbb{C}^n$ and
 $ \bigwedge(m)$ is the exterior (or Grassmann) algebra over $\mathbb C$ with $m$ generators.
 A {\it complex-analytic supermanifold} of dimension $n|m$ is a $\mathbb{Z}_2$-graded ringed space $\mathcal{M} = (\mathcal{M}_0,{\mathcal
 	O}_{\mathcal{M}})$ that
 is locally isomorphic to a complex superdomain of dimension $n|m$.
 An example of a complex-analytic supermanifold is the ringed space $(\mathcal{M}_0,\bigwedge\mathcal{E})$, where 
  $\mathcal{M}_0$ is a complex-analytic manifold and $\mathcal{E}$ is a holomorphic locally free sheaf on $\mathcal{M}_0$. In this case $\dim\,\mathcal{M} = n|m$, where $n =  \dim \mathcal{M}_0$ and $m$ is the rank of  $\mathcal{E}$.

 Let $\mathcal{M} = (\mathcal{M}_0,{\mathcal
 	O}_{\mathcal{M}})$ be a complex-analytic supermanifold
 and $ \mathcal{J}_{\mathcal{M}} = ({\mathcal
 	O}_{\mathcal{M}})_{\bar 1} + ({\mathcal
 	O}_{\mathcal{M}})_{\bar 1}^2
 $
 be the subsheaf of ideals generated
 by odd
 elements in ${\mathcal O}_{\mathcal{M}}$. We set $\mathcal{F}_{\mathcal{M}}:= {\mathcal
 	O}_{\mathcal{M}}/\mathcal{J}_{\mathcal{M}}$. Then $(\mathcal{M}_0,
 \mathcal{F}_{\mathcal{M}})$ is a usual complex-analytic manifold.  It
 is called the {\it reduction} or {\it underlying space} of $\mathcal{M}$.
 Usually we will write $\mathcal{M}_0$ instead of $(\mathcal{M}_0,
 \mathcal{F}_{\mathcal{M}})$ and $\mathcal{F}_{\mathcal{M}_0}$ instead of $\mathcal{F}_{\mathcal{M}}$. A morphism of
 supermanifolds is a morphism of the corresponding
 $\mathbb{Z}_{2}$-graded ringed spaces. If $f: \mathcal{M} \to
 \mathcal{N}$ is a morphism of supermanifolds, then we denote by
 $f_0$ the
 morphism of the underlying spaces $\mathcal{M}_0 \to
 \mathcal{N}_0$ and by $f^*$ the  morphism of the structure sheaves
 $\mathcal{O}_{\mathcal{N}} \to (f_0)_*(\mathcal{O}_{\mathcal{M}})$.

 Denote by ${\mathcal T}_{\mathcal{M}}$
 the {\it tangent sheaf} or the {\it sheaf of vector fields} of
 $\mathcal{M}$. Since the sheaf ${\mathcal O}_{\mathcal{M}}$ is
 $\mathbb{Z}_2$-graded, the tangent sheaf ${\mathcal
 	T}_{\mathcal{M}}$ is also $\mathbb{Z}_2$-graded.
Furthermore we have the following filtration in the structure sheaf $\mathcal{O}_{\mathcal{M}}$:
$$
\mathcal{O}_{\mathcal{M}}=\mathcal{J}_{\mathcal{M}}^0\supset \mathcal{J}^1_{\mathcal{M}}\supset
\mathcal{J}^2_{\mathcal{M}}\supset \ldots .
$$
This filtration induces  the following filtration in  ${\mathcal T}_{\mathcal{M}}$ 
$$
{\mathcal T}_{\mathcal{M}} = ({\mathcal T}_{\mathcal{M}})_{(-1)}\supset ({\mathcal T}_{\mathcal{M}})_{(0)} \supset ({\mathcal T}_{\mathcal{M}})_{(1)}
\supset\cdots,
$$
where
$$
({\mathcal T}_{\mathcal{M}})_{(p)}=\{ v\in {\mathcal T}_{\mathcal{M}} \mid  v(\mathcal{O}_{\mathcal{M}})\subset
\mathcal{J}^{p}_{\mathcal{M}},\, v(\mathcal{J}_{\mathcal{M}}) \subset \mathcal{J}^{p+1}_{\mathcal{M}}\} \,\,\,
\text{for}\,\,p\geq 0
$$

 Here and everywhere in the paper for simplicity we use the notation $\mathcal{U}_i:= (U_i,\mathcal{O}_{\mathcal{N}}|_{U_i})$ for a chart on a supermanifold $\mathcal N$. We say that $\{\mathcal{U}_i\}$ is an open covering of $\mathcal N$ if $\{U_i\}$ is an open covering of $\mathcal{N}_0$.

Let $\mathcal{M}$, $\mathcal{B}$ and $\mathcal{S}$ be
complex-analytic supermanifolds. 

\begin{definition}\label{de_superbundle} 
The quadruple $(\mathcal{M}, \mathcal{B}, \pi,\mathcal{S})$ is
called a {\it superbundle with total space $\mathcal{M}$, base space $\mathcal{B}$, projection $\pi: \mathcal{M} \to
	\mathcal{B}$ and with fiber $\mathcal{S}$,}  if  $\pi: \mathcal{M} \to
\mathcal{B}$ is locally trivial. In other words 
there exists an open covering $\{\mathcal{U}_i\}$ of the supermanifold $\mathcal{B}$ and
isomorphisms $\psi_i: \pi^{-1}(\mathcal{U}_i)\rightarrow
\mathcal{U}_i\times \mathcal{S}$ such that the following
diagram is commutative:
$$
\begin{CD}
\pi^{-1}(\mathcal{U}_i) @>{\psi_i}>>\mathcal{U}_i \times
\mathcal{S}\\
@V{ \pi}VV @VV {\pr_{\mathcal{U}_i}}V\\
\mathcal{U}_i@=\mathcal{U}_i
\end{CD},
$$
where $\pr_{\mathcal{U}_i}$ is the projection onto the first factor.	
\end{definition}

 Here we denote by $\pi^{-1}(\mathcal{U}_i)$ the following supermanifold $(\pi_0^{-1}(U_i), \mathcal{O}_{\mathcal{M}}|_{\pi_0^{-1}(U_i)})$. Sometimes for simplicity we will denote the superbundle $(\mathcal{M}, \mathcal{B}, \pi,\mathcal{S})$ just by $\mathcal{M}$.

Let $\mathcal{M}$ be a superbundle. 

\begin{definition}\label{de_projectable vector fields} 
A global vector field $v\in H^0(\mathcal{M}_0, \mathcal T_{\mathcal{M}})$ is called {\it projectible} with respect to $\pi$ if there exists a vector field
$v'\in H^0(\mathcal{B}_0, \mathcal T_{\mathcal{B}})$ such that
$$
v(\pi^*(f)) = v'(f) \quad\text{for all} \quad f\in \mathcal{O}_{\mathcal{B}}.
$$
The vector field $v'$ is called the {\it projection} of $v$. If
$v'= 0$, the vector field $v$ is called {\it vertical}. All vertical vector fields form a Lie subsuperalgebra in $H^0(\mathcal{M}_0, \mathcal T_{\mathcal{M}})$. 
\end{definition}

The following theorem was proven in \cite{bash}:

\begin{theorem}\label{bash}   Let $\pi: \mathcal{M}\to \mathcal{B}$
be the projection of a superbundle with fiber $\mathcal{S}$. If
$H^0(\mathcal S_0, \mathcal O_{\mathcal{S}}) \simeq \mathbb C$, then any vector field  on
$\mathcal{M}$ is projectible with respect to $\pi$.
\end{theorem}

\subsection{Split supermanifolds}\label{sec Split supermanifolds}

A supermanifold $\mathcal{M}$ is called {\it split} if
${\mathcal O}_{\mathcal{M}}\simeq \bigwedge\mathcal{E}$ for a
locally free sheaf $\mathcal{E}$ on $\mathcal{M}_0$. In this case
$\mathcal{O}_{\mathcal{M}}$ is a $\mathbb{Z}$-graded sheaf. More precisely, $\mathcal{O}_{\mathcal{M}}= \bigoplus_p (\mathcal{O}_{\mathcal{M}})_p$, where $(\mathcal{O}_{\mathcal{M}})_p$ is the image of $\bigwedge^p\mathcal{E}$ by the isomorphism $\bigwedge\mathcal{E}\to \mathcal{O}_{\mathcal{M}}$. Moreover, $\mathcal{O}_{\mathcal{M}}$ is a sheaf of $\mathcal{F}_{\mathcal M}$-modules, where as above $\mathcal{F}_{\mathcal M}$ is the structure sheaf of the underlying space $\mathcal{M}_0$. Note that $\mathcal{O}_{\mathcal{M}}\simeq \bigwedge (\mathcal{O}_{\mathcal{M}})_1$. If $\mathcal{M}$ is split, then the tangent sheaf $\mathcal{T}_{\mathcal{M}}$ is also a locally free
sheaf of $\mathcal{F}_{\mathcal M}$-modules. Furthermore, the sheaf $\mathcal{T}_{\mathcal{M}}$ is $\mathbb{Z}$-graded
\begin{equation}\label{T_p=def}
(\mathcal{T}_{\mathcal{M}})_p =\{v\in \mathcal{T}_{\mathcal{M}} \mid  v(({\mathcal
	O}_{\mathcal{M}})_q)\subset
({\mathcal O}_{\mathcal{M}})_{p+q} \,\, \,\text{for all}\,\, q\geq 0\}, \quad p\geq -1.
\end{equation}

To any split supermanifold $\mathcal M$ we can assign the following exact sequence of sheaves, see \cite{OniTransit}
\begin{equation}\label{exact sequens for T_p}
\begin{split}
&0\to (\mathcal A_{\mathcal M})_p \stackrel{\alpha}{\longrightarrow}
(\mathcal{T}_{\mathcal{M}})_p \stackrel{\beta}{\longrightarrow} (\mathcal C_{\mathcal M})_p\to 0,
\end{split}
\end{equation}
where
$$
(\mathcal A_{\mathcal M})_p:= ({\mathcal O}_{\mathcal{M}})_1^*\otimes \bigwedge^{p+1}({\mathcal O}_{\mathcal{M}})_1\quad \text{and} \quad (\mathcal C_{\mathcal M})_p:= \Theta_{\mathcal{M}_0}\otimes \bigwedge^{p}({\mathcal O}_{\mathcal{M}})_1.
$$
Here $\Theta_{\mathcal{M}_0}$ is the tangent sheaf of the manifold $\mathcal{M}_0$, the map $\alpha$ assigns the corresponding vector field from $(\mathcal{T}_{\mathcal{M}})_p$ to any morphism $({\mathcal O}_{\mathcal{M}})_1 \to  \bigwedge^{p+1}({\mathcal O}_{\mathcal{M}})_1 $, and $\beta$ is the restriction of a vector field $v\in (\mathcal{T}_{\mathcal{M}})_p$ to the subsheaf $\mathcal{F}_{\mathcal M}$.

Split supermanifolds form a category. More precisely, objects of this category are all split supermanifolds $\mathcal{M}$ with a fixed
isomorphism $\mathcal{O}_{\mathcal{M}}\simeq \bigwedge \mathcal{E}$ for a certain
locally free sheaf $\mathcal{E}$ on $\mathcal{M}_0$, and morphisms are all morphisms of ringed spaces preserving the
$\mathbb{Z}$-gradings. There is a functor $\gr$ from the category of supermanifolds to the
category of split supermanifolds. Let us briefly describe this
 construction. Let $\mathcal{M}$ be any supermanifold.
As above, denote by $\mathcal{J}_{\mathcal{M}}\subset \mathcal{O}_{\mathcal{M}}$ the subsheaf of
ideals generated by odd elements of $\mathcal{O}_{\mathcal{M}}$. Then by
definition $\gr\mathcal{M}$ is the split supermanifold with the structure sheaf
$$
\mathcal{O}_{\gr\mathcal{M}}= \bigoplus_{p\geq 0} (\mathcal{O}_{\gr\mathcal{M}})_p,
\quad (\mathcal{O}_{\gr\mathcal{M}})_p:= \mathcal{J}_{\mathcal{M}}^p/\mathcal{J}_{\mathcal{M}}^{p+1},\quad
\mathcal{J}_{\mathcal{M}}^0:=\mathcal{O}_{\mathcal{M}}.
$$
In this case $(\mathcal{O}_{\gr\mathcal{M}})_1$ is a locally free sheaf and
there is a natural isomorphism of $\mathcal{O}_{\gr\mathcal{M}}$ onto
$\bigwedge (\mathcal{O}_{\gr\mathcal{M}})_1$.
If $\psi=(\psi_{0},\psi^*):\mathcal{M}\to
\mathcal{N}$ is a morphism of supermanifolds, then
$\gr(\psi)=(\psi_{0},\gr(\psi^*))$, where $\gr(\psi^*):
\mathcal{O}_{\mathcal{N}} \to \mathcal{O}_{\mathcal{M}}$ is defined by
$$
\gr(\psi^*)(f+\mathcal{J}_{\mathcal{N}}^p): = \psi^*(f)+\mathcal{J}_{\mathcal{M}}^p
\,\,\text{for}\,\, f\in (\mathcal{J}_{\mathcal{N}})^{p-1}.
$$
Recall that by definition every morphism $\psi$ of supermanifolds is
even and as a consequence $\psi*$ sends $\mathcal{J}_{\mathcal{N}}^p$ into
$\mathcal{J}_{\mathcal{M}}^p$.

If $\mathcal{T}_{\mathcal{M}}$ is the tangent sheaf of a supermanifold
$\mathcal{M}$ and ${\mathcal{T}}_{\gr\mathcal{M}}$ is the tangent sheaf of
$\gr \mathcal{M}$, then as above the tangent sheaf ${\mathcal{T}}_{\gr\mathcal{M}}$ is $\mathbb Z$-graded
$$
{\mathcal{T}}_{\gr\mathcal{M}}=\bigoplus_{p\geq -1}({\mathcal{T}}_{\gr\mathcal{M}})_{p},
$$
and
$$
({\mathcal{T}}_{\gr\mathcal{M}})_p\simeq (\mathcal{T}_{\mathcal{M}})_{(p)}/(\mathcal{T}_{\mathcal{M}})_{(p+1)}
\,\,\,\text{for $p\geq -1$, see \cite{Oni_contruc}.}
$$

A {\it split superbundle} is a superbundle in the category of split
supermanifolds. More precisely, it is a superbundle $(\mathcal{M}, \mathcal{B}, \pi, \mathcal{S})$, where $\mathcal{M}$, $\mathcal{B}$ and $\mathcal{S}$ are
split and  $\pi$ preserves the
$\mathbb{Z}$-gradings. If $\pi:\mathcal{M}\to \mathcal{B}$
is a superbundle with fiber $\mathcal{S}$ then
$\gr\pi:\gr\mathcal{M}\to \gr\mathcal{B}$ is a split
superbundle with fiber $\gr\mathcal{S}$.

\subsection{Retract of a split superbundle and its tangent sheaf}\label{sec Retract of a split superbundle}

In this section $(\mathcal{M},\mathcal{B},\pi,\mathcal{S})$ is a split superbundle, this is $\mathcal{M}$, $\mathcal{B}$ and $\mathcal{S}$ are split supermanifolds with fixed $\mathbb Z$-gradings and $\pi^*: \mathcal{O}_{\mathcal{B}}\to (\pi_0)_*(\mathcal{O}_{\mathcal{M}})$ preserves these $\mathbb Z$-gradings. To simplify notations we write in this section $\mathcal{O}$ instead of $\mathcal{O}_{\mathcal{M}}$ and $\mathcal{T}$ instead of $\mathcal{T}_{\mathcal{M}}$. 

 Denote by $\mathcal{J}_{\mathcal{B}}$ the sheaf
of ideals in $\mathcal{O}$ generated by 
$\pi^*((\mathcal{O}_{\mathcal{B}})_{\bar 1})$. By definition of a split superbundle
$\mathcal{J}_{\mathcal{B}}$ is a sheaf of $\mathbb{Z}$-graded ideals. We have the following filtration in $\mathcal{O}$
\begin{equation}\label{filtr_J_B}
\mathcal O=\mathcal J_{\mathcal{B}}^0\supset \mathcal J_{\mathcal{B}}^1\supset \mathcal
J_{\mathcal{B}}^2\supset \dots .
\end{equation}
 Since $\mathcal{M}$ is split, we have $\mathcal{O}=\bigoplus\limits_p \mathcal{O}_p$ and Filtration
(\ref{filtr_J_B}) gives rise to the following filtration in any $\mathcal O_p$
\begin{equation}
\label{form_filtr O_(p(q))}
\begin{split}
\mathcal O_p=\mathcal O_{p(0)}\supset \mathcal O_{p(1)} \supset\dots
\supset\mathcal O_{p(p)}\supset \mathcal O_{p(p+1)}=0,
\end{split}
\end{equation}
where $\mathcal O_{p(q)}= \mathcal J^q_{\mathcal{B}}\cap \mathcal O_p$. 
We set
$$
\widehat{\mathcal O}_{pq}:=\mathcal O_{p(q)}/ \mathcal O_{p(q+1)} \quad\text{and} \quad \widehat{\mathcal O}:= \bigoplus\limits_{p,q\geq0}\widehat{\mathcal
	O}_{pq}.
$$
The sheaf $\widehat{\mathcal O}$ is a sheaf of superalgebras with respect to the following multiplication: 
$$
(f+ \mathcal O_{p(q+1)}) (g+ \mathcal O_{p'(q'+1)}) = fg + \mathcal O_{p+p',(q+q'+1)},
$$
where $f\in  \mathcal O_{p(q)}$ and $g\in  \mathcal O_{p'(q')}$. The sheaves $\widehat{\mathcal O}_{11}$ and $\widehat{\mathcal O}_{10}$ are locally free. Moreover, 
$$
\widehat{\mathcal O}_{pq}= \bigwedge^{p-q}
\widehat{\mathcal O}_{10}\otimes \bigwedge^{q} \widehat{\mathcal O}_{11}.
$$ 
In particular the ringed space $\widehat{\mathcal M}:= (\mathcal M_0,\widehat{\mathcal O})$ is a split supermanifold. 
Since, $\pi^*((\mathcal{O}_{\mathcal B})_p) \subset \mathcal O_{p(p)}= \widehat{\mathcal O}_{pp}$, the morphism
$\widehat{\pi}: \widehat{\mathcal M}\to \mathcal{B}$ is defined, where 
$\widehat{\pi}^* (f) = \pi^*(f)$ for $f\in \mathcal{O}_{\mathcal{B}}$. The morphism $\widehat{\pi}: \widehat{\mathcal M}\to
\mathcal{B}$ preserves the fixed $\mathbb Z$-grading in $\mathcal{O}_{\mathcal{B}}$ and the $\mathbb Z$-grading $ \widehat{\mathcal O}_p:=\bigoplus\limits_{q}\widehat{\mathcal
	O}_{pq}$ in $\widehat{\mathcal
	O}$. Hence, $\widehat{\mathcal M}$ is a split superbundle. We will call the superbundle $\widehat{\pi}: \widehat{\mathcal M}\to
\mathcal{B}$ the {\it
retract of the split superbundle} $\pi: \mathcal{M}\to
\mathcal{B}$.

Further consider the tangent sheaf $\mathcal{T}$ of the split supermanifold $\mathcal M$. The filtration
(\ref{filtr_J_B}) induces the following filtration in
$\mathcal{T}_p$
\begin{equation}\label{eq_filtration in T_p}
\mathcal T_p=\mathcal T_{p(-1)}\supset \mathcal T_{p(0)}\supset
\dots\mathcal T_{p(q)}\supset\dots\mathcal T_{p(p+1)}\supset\mathcal
T_{p(p+2)}=0.
\end{equation}
where
$$
\mathcal T_{p(q)}= \lbrace v\in \mathcal T_p\,\,\vert\,\, v(\mathcal J_{\mathcal B})\subset \mathcal J_{\mathcal B}^{q+1}, \,\,\, v(\mathcal O)\subset \mathcal J_{\mathcal B}^{q} \rbrace, \quad q\geq -1.
$$
We put
$$
\widehat{\mathcal T}= \bigoplus_{p,q\geq -1}\widehat{\mathcal
T}_{pq},\,\, \text{where}\,\,\widehat{\mathcal T}_{pq}=\mathcal
T_{p(q)}/\mathcal T_{p(q+1)}.
$$
The sheaf $\widehat{\mathcal T}$ possesses  a Lie superalgebra structure that is induced by the Lie superalgebra structure on ${\mathcal T}$. Moreover,  $\widehat{\mathcal T}$ has a natural geometric interpretation. Denote by $Der \widehat{\mathcal O}$ the sheaf of derivations of
$\widehat{\mathcal O}$. It is a sheaf of Lie superalgebras. Clearly,
$$
Der \widehat{\mathcal O}= \bigoplus_{p,q\geq -1} Der_{pq}
\widehat{\mathcal O},
$$
 where
$$
Der_{pq} \widehat{\mathcal O}= \{ u\in Der \widehat{\mathcal O} \,\,
\mid\,\, u(\widehat{\mathcal O}_{st})\subset \widehat{\mathcal
O}_{s+p, t+q}, \forall s,t\in \mathbb Z \}, \,\,\,p,q\geq -1,
$$

\begin{lemma}\label{Lemma_T isom Der}
	We have $\widehat{\mathcal
T}\simeq Der  \widehat{\mathcal O}$ as sheaves of Lie superalgebras and $\widehat{\mathcal T}_{pq}\simeq Der_{pq}\widehat{\mathcal O}$.
\end{lemma}

\medskip

\noindent{Proof.} Let us take $u\in (\mathcal T_{p(q)})_x$, where $x\in
\mathcal M_0$. Since 
$$
u((\mathcal O_{r(s)})_x)\subset (\mathcal
O_{r+p,(s+q)})_x \quad  \text{and} \quad u((\mathcal O_{r(s+1)})_x)\subset (\mathcal
O_{r+p,(s+q+1)})_x,
$$
we see that $u$ determines the $(p,q)$-derivation
$\widehat{u}: (\widehat{\mathcal O}_{rs})_x\to (\widehat{\mathcal
O}_{r+p,s+q})_x$. Further, it is easy to see that $\widehat{u}=0$ if and only if $u\in (\mathcal
T_{p(q+1)})_x$. Therefore, the map 
$$
\alpha_{pq}:\widehat{\mathcal
T}_{pq}\rightarrow Der_{pq} \widehat{\mathcal O},\quad u+(\mathcal
T_{p(q+1)})_x\mapsto \widehat{u}
$$
is injective. Clearly, locally we can find pre-images of basic
elements. Hence, $\alpha_{pq}$ is surjective. Moreover, by definitions of all structures, we see that
$$
\bigoplus_{pq}\alpha_{pq}: \bigoplus_{pq} \widehat{\mathcal T}_{pq}\to
\bigoplus_{pq}Der_{pq}\widehat{\mathcal O}
$$
 is a homomorphism of
sheaves of Lie superalgebras.$\Box$

\medskip

We will identify $\widehat{\mathcal T}$ and $Der\widehat{\mathcal
O}$ using the isomorphism $\oplus \alpha_{pq}$ from Lemma \ref{Lemma_T isom Der}.

Let us rewrite 
exact Sequence (\ref{exact sequens for T_p})  for $\widehat{\mathcal T}_p$. Recall that $ \widehat{\mathcal O}_p=\bigoplus\limits_{q}\widehat{\mathcal
	O}_{pq}$.  We get
$$
0\to (\widehat{\mathcal O}_{10}\oplus\widehat{\mathcal
O}_{11})^*\otimes\bigwedge^{p+1}(\widehat{\mathcal O}_{10}
\oplus\widehat{\mathcal O}_{11})\to \widehat{\mathcal T}_p\to
\Theta\otimes \bigwedge^{p} (\widehat{\mathcal
O}_{10}\oplus\widehat{\mathcal O}_{11})\to 0.
$$
Here $\Theta$ is the tangent sheaf of the manifold $\mathcal{M}_0$
This exact sequence is a direct sum of the following exact sequences.
\begin{equation}\label{posled dlja T_pq}
0\to \mathcal{A}_{pq}\to \widehat{\mathcal T}_{pq}\to
\mathcal{C}_{pq}\to 0,
\end{equation}
where
\begin{align*}
\mathcal{A}_{pq}&:= \widehat{\mathcal{O}}_{10}^*\otimes\bigwedge^{p-
	q+1}\widehat{\mathcal{O}}_{10}\otimes
\bigwedge^q\widehat{\mathcal{O}}_{11} + \widehat{\mathcal{O}}_{11}^*\otimes\bigwedge^{p- q}\widehat{\mathcal{O}}_{10}\otimes
\bigwedge^{q+1}\widehat{\mathcal{O}}_{11},\\
\mathcal{C}_{pq}&:= \Theta\otimes\bigwedge^{p-q}\widehat{\mathcal{O}}_{10}\otimes \bigwedge^{q}\widehat{\mathcal
	O}_{11}.
\end{align*}

The sheaves $\mathcal{A}_{pq}$ and $\mathcal{C}_{pq}$ possess another description in terms of factor sheaves. More precisely, by definition we have 
$$
\mathcal A_p= \mathcal{O}_1^*\otimes \bigwedge^{p+1} \mathcal{O}_1 \quad \text{and} \quad \mathcal C_p= \Theta\otimes \bigwedge^{p} \mathcal{O}_1,
$$
see (\ref{exact sequens for T_p}).
These sheaves possess the following filtrations
$$
\mathcal A_p = \mathcal A_{p(-1)}\supset \mathcal A_{p(0)} \supset\ldots,\quad \quad
\mathcal C_p = \mathcal C_{p(0)}\supset \mathcal C_{p(1)} \supset\ldots,\\
$$
where
$$
\mathcal A_{p(q)}=\mathcal{O}^*_{10}\otimes \bigwedge^{p-q+1} \mathcal
O_{1(0)}\otimes \bigwedge^{q}\mathcal O_{1(1)} + \mathcal{O}_1^*\otimes
\bigwedge^{p-q}\mathcal O_{1(0)}\otimes \bigwedge^{q+1}\mathcal
O_{1(1)} 
$$
 and
$$
\mathcal C_{p(q)}=\Theta\otimes \bigwedge^{p-q}\mathcal
O_{1(0)}\otimes \bigwedge^{q}\mathcal O_{1(1)}.
$$
Note that $\mathcal{O}^*_{10}\subset \mathcal{O}_1^*$.

\begin{lemma}\label{Lemma_A_pq= B-pq=}  
	We have $\mathcal A_{pq}=\mathcal
A_{p(q)}/\mathcal A_{p(q+1)}$ and $\mathcal C_{pq}=\mathcal
C_{p(q)}/\mathcal C_{p(q+1)}$.
\end{lemma}

\medskip

\noindent{Proof.} The idea of the proof is similar to the idea of the proof of Lemma \ref{Lemma_T isom Der}.$\Box$

\section{ Exact sequences associated to the tangent sheaf of a split superbundle}

Let $\pi:\mathcal{M}\to \mathcal{B}$ be an arbitrary
superbundle with fiber $\mathcal{S}$ and
$\gr\pi:\gr \mathcal{M}\to \gr\mathcal{B}$ be
the corresponding split superbundle with fiber $\gr\mathcal{S}$.
Denote also by $\widehat{\pi}:\widehat{\mathcal{M}}\to
\widehat{\mathcal{B}}$ the retract of
$\gr\pi:\gr \mathcal{M}\to \gr\mathcal{B}$. 
From now on we fix the following notations
\begin{equation}\label{def notations O_M, O_B, O_S}
\begin{split}
\mathcal O,\quad \tilde{\mathcal O}= \bigoplus_{p\geq 0} \tilde{\mathcal O}_p, \quad \widehat{\mathcal O} = \bigoplus_{p,q\geq 0} \widehat{\mathcal O}_{pq}, \quad \mathcal O_{\mathcal B},\quad \tilde{\mathcal O}_{\mathcal B}=\widehat{\mathcal O}_{\mathcal B}= \bigoplus_{p\geq 0} (\tilde{\mathcal O}_{\mathcal B})_p, \\
\mathcal O_{\mathcal S},\quad \tilde{\mathcal O}_{\mathcal S}=\widehat{\mathcal O}_{\mathcal S}= \bigoplus_{p\geq 0} (\tilde{\mathcal O}_{\mathcal S})_p
\end{split}
\end{equation}
for the structure sheaves of $\mathcal M$, $\gr \mathcal{M}$, $\widehat{\mathcal{M}}$, $\mathcal B$, $\gr \mathcal{B} =\widehat{\mathcal{B}}$, $\mathcal S$ and $\gr \mathcal{S} =\widehat{\mathcal{S}}$, respectively. To simplify notations sometimes we will denote the bundle projections $\pi$, $\gr{\pi}$ and $\widehat{\pi}$ simply by $\pi$. From now on we also denote by
\begin{equation}\label{def notations T_M, T_B, T_S}
\begin{split}
\mathcal T,\quad  \tilde{\mathcal T} = \bigoplus_{p\geq -1}\tilde{\mathcal T}_p,\quad \text{and} \quad \widehat{\mathcal T} =
\bigoplus_{p,q \geq -1}\widehat{\mathcal T}_{pq}, \quad \mathcal T_{\mathcal B},\quad \tilde{\mathcal T}_{\mathcal B}=\widehat{\mathcal T}_{\mathcal B}= \bigoplus_{p\geq -1} (\tilde{\mathcal T}_{\mathcal B})_p, \\
\mathcal T_{\mathcal S},\quad \tilde{\mathcal T}_{\mathcal S}=\widehat{\mathcal T}_{\mathcal S}= \bigoplus_{p\geq -1} (\tilde{\mathcal T}_{\mathcal S})_p
\end{split}
\end{equation}
the tangent
sheaves of $\mathcal M$, $\gr \mathcal{M}$, $\widehat{\mathcal{M}}$, $\mathcal B$, $\gr \mathcal{B} =\widehat{\mathcal{B}}$, $\mathcal S$ and $\gr \mathcal{S} =\widehat{\mathcal{S}}$, respectively. A local vector field $v\in \mathcal T$ is called {\it vertical} if $v(\pi^*(\mathcal{O}_{\mathcal{B}})) = 0$. Clearly, the commutator of two vertical vector fields is again a  vertical vector field.
 Denote by $\mathcal
T^v\subset \mathcal T$, $\tilde{\mathcal T}_p^v\subset
\tilde{\mathcal T}_p$ and $\widehat{\mathcal T}_{pq}^v\subset
\widehat{\mathcal T}_{pq}$ the subsheaves of vertical vector fields, respectively. The sheaves $\mathcal T^v$ and $\tilde{\mathcal T}_p^v$ possess the following filtrations:
$$
\mathcal T^v = \mathcal T^v_{(-1)}\supset \mathcal T^v_{(0)}\supset
\cdots,\quad \quad
\tilde{\mathcal T}^v_p = \tilde{\mathcal T}^v_{p(-1)}\supset
\tilde{\mathcal T}^v_{p(0)}\supset \cdots,
$$
where $\mathcal T^v_{(p)}:= \mathcal T^v\cap \mathcal T_{(p)}$ and
$\tilde{\mathcal T}^v_{p(q)}:= \tilde{\mathcal T}_p^v\cap \mathcal T_{p(q)}$. 

\begin{lemma}\label{Lemma_T^v_p= T^v_pq=}  We have $\tilde{\mathcal T}_p^v \simeq \mathcal T_{(p)}^v/  \mathcal T_{(p+1)}^v$ and $\widehat{\mathcal T}_{pq}^v \simeq \tilde{\mathcal T}_{p(q)}^v/ \tilde{\mathcal T}_{p(q+1)}^v$.
\end{lemma}

\medskip

\noindent{Proof.} Consider the natural map $\mathcal
T_{(p)}^v \to \tilde{\mathcal
T}_p^v$. It is surjective and its kernel is $\mathcal T_{(p+1)}^v$.
For $\widehat{\mathcal T}_{pq}^v$ the argument is similar.$\Box$

\medskip

We put 
$$
\mathcal T^h= \mathcal T/ \mathcal T^v,\quad \tilde{\mathcal
T}^h_{p}:= \tilde{\mathcal T}_{p}/
 \tilde{\mathcal
T}^v_{p}\quad \text{and}\quad \widehat{\mathcal T}^h_{pq}:= \widehat{\mathcal
T}_{pq}/ \widehat{\mathcal T}^v_{pq}.
 $$
 Denote also 
 $$
 \mathcal T^h_{(p)}= \mathcal T_{(p)}/ \mathcal
 T^v_{(p)}\quad  \text{and}\quad \tilde{\mathcal T}^h_{p(q)}= \tilde{\mathcal T}_{p(q)}/ \tilde{\mathcal T}^v_{p(q)}.
 $$
Since we have the natural inclusions $\mathcal T_{(p+1)} \hookrightarrow \mathcal T_{(p)}$ and  $\mathcal T^v_{(p+1)} \hookrightarrow \mathcal T^v_{(p)}$, we can define the sheaf morphism $\Phi:\mathcal T^h_{(p+1)} \to  \mathcal T^h_{(p)}$. Furthermore, the sheaf morphisms $\mathcal T_{(p)} \to \tilde{\mathcal T}_{p}$ and $\mathcal T_{(p)}^v \to \tilde{\mathcal T}^v_{p}$ give rise to the morphism $\Psi: \mathcal T_{(p)}^h \to \tilde{\mathcal T}^h_{p}$. 
We will need the following lemmas:

\begin{lemma}\label{Lemma_diagram tilde T^h commutative}  The following diagram is commutative,  lines and 
 columns are exact:
  \begin{center}
  \begin{tikzpicture}
 \matrix (m) [matrix of math nodes,row sep=1em,column sep=1em,minimum width=1em]
 {& 0& 0&0 & \\
 	 0 &	\tilde{\mathcal T}^v_{p} &  \tilde{\mathcal T}_{p} & \tilde{\mathcal T}_{p}^h & 0\\
 0  &	\stackrel{}{\mathcal T^v_{(p)}} & \mathcal T_{(p)} & \mathcal T_{(p)}^h & 0 \\
 0 &	\stackrel{}{\mathcal T^v_{(p+1)}} & \mathcal T_{(p+1)} & \mathcal T_{(p+1)}^h & 0 \\
 & 0& 0&0 & \\};
 \draw[->]
  (m-2-1) edge  (m-2-2)
  (m-3-1) edge  (m-3-2)
  (m-4-1) edge  (m-4-2)
   (m-2-2) edge node [left] {} (m-1-2)
 edge  node [below] {$\chi$} (m-2-3)
 (m-3-2) edge node [right] {$\sigma$} (m-2-2)
         edge  node [below] {$\theta$} (m-3-3)
 (m-4-2) edge node [left] {} (m-3-2)
         edge  node [below] {} (m-4-3)
  (m-5-2) edge  (m-4-2)  
 (m-5-3) edge  (m-4-3)
 (m-5-4) edge  (m-4-4) 
 
 (m-2-3) edge  (m-1-3)
         edge  node [below] {$\delta$} (m-2-4)
       
 (m-3-3) edge  node [right] {$\epsilon$} (m-2-3)
         edge  node [below] {$\gamma$} (m-3-4)
  (m-4-3) edge  (m-3-3)
          edge  node [below] {$\lambda$} (m-4-4)
          
 (m-2-4) edge (m-1-4)
         edge (m-2-5)
 (m-3-4) edge node [right] {$\Psi$} (m-2-4)
         edge (m-3-5)
 (m-4-4) edge  node [right] {$\Phi$} (m-3-4)
         edge   (m-4-5);
 \end{tikzpicture}
 \end{center}
\end{lemma}

\noindent{Proof.} The exactness of all lines and first two columns follows from definitions and Lemma \ref{Lemma_T^v_p= T^v_pq=}. Let us show the exactness of
  $$
  0\to \mathcal T_{(p+1)}^h \stackrel{\Phi}{\longrightarrow} \mathcal T_{(p)}^h \stackrel{\Psi}{\longrightarrow} \tilde{\mathcal T}_{p}^h \to 0
  $$
 for example in the term $\mathcal T_{(p)}^h$. Let us take $v\in \mathcal T_{(p)}^h$ such that $\Psi (v)= 0$. We have to show that there exists $w\in \mathcal T_{(p+1)}^h$ such that $\Phi(w)=v$. Using commutativity and exactness we get the following. There exists $x\in \mathcal T_{(p)}$ such that $\gamma(x)=v$. Since $\delta\circ \epsilon (x)=0$, there exists $z\in \tilde{\mathcal T}_{p}^v$ such that $\chi(z)= \epsilon(x)$. Further there exists $t\in \mathcal T_{(p)}^v$ such that $\sigma(t)= z$. Consider $\tilde x= x-\theta(t)$. We have $\gamma(\tilde x) = \gamma(x)- \gamma(\theta(t))= \gamma(x)=v$. Hence, $\tilde x$ is also a representant of $v$ in $\mathcal T_{(p)}$. Moreover, $\epsilon(\tilde x)= \epsilon(x)- \epsilon(\theta(t)) = 0$, therefore, $\exists t\in \mathcal T_{(p+1)}$ such that we can set $w:=\lambda(r)$. 
   This observation completes the proof.$\Box$

\begin{lemma}\label{Lemma_diagram tilde T^h_pq commutative} The following diagram is commutative,  lines and 
  	columns are exact:
  \begin{center}
  	\begin{tikzpicture}
  	\matrix (m) [matrix of math nodes,row sep=1em,column sep=1em,minimum width=1em]
  	{& 0& 0&0 & \\
  		0 &	\widehat{\mathcal T}^v_{pq} &  \widehat{\mathcal T}_{pq} & \widehat{\mathcal T}_{pq}^h & 0\\
  		0  &	\tilde{\mathcal T}^v_{p(q)} & \tilde{\mathcal T}_{p(q)} & \tilde{\mathcal T}_{p(q)}^h & 0 \\
  		0 &	\tilde{\mathcal T}^v_{p(q+1)} & \tilde{\mathcal T}_{p(q+1)} & \tilde{\mathcal T}^h_{p(q+1)} & 0 \\
  		& 0& 0&0 & \\};
  	\draw[->]
  	(m-2-1) edge  (m-2-2)
  	(m-3-1) edge  (m-3-2)
  	(m-4-1) edge  (m-4-2)
  	(m-2-2) edge  (m-1-2)
  	edge   (m-2-3)
  	(m-3-2) edge  (m-2-2)
  	edge   (m-3-3)
  	(m-4-2) edge  (m-3-2)
  	edge  (m-4-3)
  	(m-5-2) edge  (m-4-2)  
  	(m-5-3) edge  (m-4-3)
  	(m-5-4) edge  (m-4-4) 
  	
  	(m-2-3) edge  (m-1-3)
  	edge   (m-2-4)
  	
  	(m-3-3) edge  (m-2-3)
  	edge  (m-3-4)
  	(m-4-3) edge  (m-3-3)
  	edge  (m-4-4)
  	
  	(m-2-4) edge (m-1-4)
  	edge (m-2-5)
  	(m-3-4) edge  (m-2-4)
  	edge (m-3-5)
  	(m-4-4) edge   (m-3-4)
  	edge   (m-4-5);
  	\end{tikzpicture}
  \end{center}
  
\end{lemma}

\medskip

\noindent{Proof.} The proof is similar to the proof of Lemma \ref{Lemma_diagram tilde T^h commutative}.$\Box$

We put 
$$
\mathcal A_{pq}^v= \widehat{\mathcal{O}}_{10}^*\otimes
\bigwedge^q\widehat{\mathcal{O}}_{11}\otimes\bigwedge^{p-
q+1}\widehat{\mathcal{O}}_{10}, \quad\quad \mathcal C_{pq}^v= \Theta^v\otimes
\bigwedge^{q}\widehat{\mathcal
O}_{11}\otimes\bigwedge^{p-q}\widehat{\mathcal{O}}_{10},
$$
 where
$\Theta^v$ is the sheaf of vertical vector fields on $M$ with respect to $\pi_{0}$.

\begin{lemma}\label{Lemma_diag for A,B}  The following diagram is
commutative, the lines and the
 columns are exact:
 \begin{center}
 	\begin{tikzpicture}
 	\matrix (m) [matrix of math nodes,row sep=1.5em,column sep=1.5em,minimum width=2em]
 	{& 0& 0&0 & \\
 	0  &	{\mathcal C}^v_{pq} & 	{\mathcal C}_{pq} & 	{\mathcal C}^h_{pq} & 0 \\
 	
 		0 &	\widehat{\mathcal T}^v_{pq} &  \widehat{\mathcal T}_{pq} & \widehat{\mathcal T}_{pq}^h & 0\\
 		
 		0  &	{\mathcal A}^v_{pq} & 	{\mathcal A}_{pq} & 	{\mathcal A}^h_{pq} &0 \\
 		& 0& 0&0 & \\};
 	\draw[->]
 	(m-2-1) edge  (m-2-2)
 	(m-3-1) edge  (m-3-2)
 	(m-4-1) edge  (m-4-2)
 	(m-2-2) edge  (m-1-2)
 	edge   (m-2-3)
 	(m-3-2) edge  (m-2-2)
 	edge   (m-3-3)
 	(m-4-2) edge  (m-3-2)
 	edge  (m-4-3)
 	(m-5-2) edge  (m-4-2)  
 	(m-5-3) edge  (m-4-3)
 	(m-5-4) edge  (m-4-4) 
 	
 	(m-2-3) edge  (m-1-3)
 	edge   (m-2-4)
 	
 	(m-3-3) edge  (m-2-3)
 	edge  (m-3-4)
 	(m-4-3) edge  (m-3-3)
 	edge  (m-4-4)
 	
 	(m-2-4) edge (m-1-4)
 	edge (m-2-5)
 	(m-3-4) edge  (m-2-4)
 	edge (m-3-5)
 	(m-4-4) edge   (m-3-4)
 	edge   (m-4-5);
 	\end{tikzpicture}
 \end{center}
 where
$$
\begin{array}{l}
\mathcal A_{pq}^h= \widehat{\mathcal{O}}_{11}^*\otimes
\bigwedge^{q+1}\widehat{\mathcal{O}}_{11}\otimes\bigwedge^{p-
q}\widehat{\mathcal{O}}_{10},\,\,\,\, \mathcal C_{pq}^h=\Theta/
\Theta^v \otimes \bigwedge^{q}\widehat{\mathcal
O}_{11}\otimes\bigwedge^{p-q}\widehat{\mathcal{O}}_{10}
\end{array}.
$$

\end{lemma}

\medskip

\noindent{Proof.} By definition, $\mathcal A_{pq}=\mathcal A_{pq}^h+ \mathcal A_{pq}^v$. Further, from exactness of the sequence
$$
0\to \Theta^v \longrightarrow \Theta \longrightarrow \Theta^h\to 0
$$
it follows the exactness of
$$
0\to \mathcal C_{pq}^v \longrightarrow \mathcal C_{pq} \longrightarrow \mathcal C_{pq}^h\to 0.
$$
The proof of exactness of
$$
0\to \mathcal A_{pq}^v \longrightarrow \widehat{\mathcal T}_{pq}^v \longrightarrow \mathcal C_{pq}^v\to 0
$$
is similar to the proof of exactness of (\ref{exact sequens for T_p}) in \cite{OniTransit}. The rest of the proof is similar to the proof of Lemma \ref{Lemma_diagram tilde T^h commutative}.$\Box$

\medskip

 We set
\begin{equation}\label{eq definition A_B,A_S,C_B,C_S} 
\begin{array}{l}
(\mathcal{A}_{\mathcal{B}})_p:= (\widehat{\mathcal{O}}_{\mathcal{B}})_1^*\otimes \bigwedge\limits^{p+1} (\widehat{\mathcal{O}}_{\mathcal{B}})_1, 
\quad (\mathcal{C}_{\mathcal{B}})_p:= \Theta_{\mathcal{B}_0}\otimes \bigwedge\limits^{p} (\widehat{\mathcal{O}}_{\mathcal{B}})_1,\\
(\mathcal{A}_{\mathcal{S}})_p:= (\widehat{\mathcal{O}}_{\mathcal{S}})_1^*\otimes
\bigwedge\limits^{p+1} (\widehat{\mathcal{O}}_{\mathcal{S}})_1,  \quad
 (\mathcal{C}_{\mathcal{S}})_p:= \Theta_{\mathcal{S}_0}\otimes \bigwedge\limits^{p} (\widehat{\mathcal{O}}_{\mathcal{S}})_1,
\end{array}
\end{equation}
where $\Theta_{\mathcal{B}_0}$ and $\Theta_{\mathcal{S}_0}$ are tangent sheaves of the manifolds $\mathcal{B}_0$ and $\mathcal{S}_0$,
respectively.  Let us choose a trivialization domain
$\mathcal{U}: = (\mathcal{U}_0 ,\tilde{\mathcal{O}}_{\mathcal{B}} |_{\mathcal{U}_0})$ of the bundle $\gr\mathcal{M}$ and denote
by $\tilde{\psi}_{\mathcal{U}}$ a trivialization isomorphism
$$
\tilde{\psi}_{\mathcal{U}}:\gr\pi^{-1}(\mathcal{U}) \rightarrow
\mathcal{U}\times \gr\mathcal{S}.
$$
For the bundle $\widehat{\mathcal{M}}$ we have the corresponding trivialization domain ${\mathcal{U}}$  of the bundle $\widehat{\mathcal{M}}$ and the following trivialization isomorphism:
$$
\widehat{\psi}_{{\mathcal{U}}}:\widehat\pi^{-1}(\mathcal{U}) \rightarrow
{\mathcal{U}}\times \widehat{\mathcal{S}}.
$$
 We will identify $\widehat\pi^{-1}(\mathcal{U})$ with $\mathcal{U}\times \widehat{\mathcal{S}}$ using this isomorphism. Denote by $\pi_{\mathcal S}$ the projection $\mathcal{U}\times \widehat{\mathcal{S}} \to \widehat{\mathcal{S}}$. 
 
A (local) description of the sheaves introduced above is given in the following proposition. 

\begin{proposition}\label{prop_geometr A,B,T}  
Let $\pi=(\pi_0,\pi^*):\mathcal M\to \mathcal B$ be a superbundle and $\widehat\pi=(\widehat\pi_0,\widehat\pi^*):\widehat{\mathcal M}\to \widehat{\mathcal B}$ be the retract of $\gr \mathcal M$. Then we have
 \begin{enumerate}
 	\item $\widehat{\mathcal O}_{11} \simeq  \widehat\pi^* ((\widehat{\mathcal O}_{\mathcal B})_1)$;
\item $\mathcal A_{pq}^h\simeq \bigwedge\limits^{p- q}\widehat{\mathcal{O}}_{10} \otimes {\pi}^*((\mathcal A_{\mathcal{B}})_{q})$ and
 $\mathcal A_{pq}^v\vert_{\widehat\pi^{-1}(\mathcal{U})} \simeq \bigwedge\limits^{q}\widehat{\mathcal O}_{11}\vert_{\mathcal{U}} \otimes \pi_{\mathcal S}^* (\mathcal
 A_{\mathcal{S}})_{p-q}$;
\item $\mathcal C_{pq}^h\simeq  \bigwedge\limits^{p- q}\widehat{\mathcal{O}}_{10}\otimes {\pi}^*((\mathcal C_{\mathcal{B}})_{q})$ and
 $\mathcal C_{pq}^v\vert_{\widehat\pi^{-1}(\mathcal{U})} \simeq \bigwedge\limits^{q}\widehat{\mathcal O}_{11}\vert_{\mathcal{U}}\otimes \pi_{\mathcal S}^* (\mathcal
 C_{\mathcal{S}})_{p-q}$;
\item $\widehat{\mathcal T}^h_{pq} \simeq \bigwedge\limits^{p- q}\widehat{\mathcal{O}}_{10}\otimes {\pi}^*((\widehat{\mathcal T}_{\mathcal{B}})_{q})$ and
$\widehat{\mathcal T}^v_{pq} \vert_{\widehat\pi^{-1}(\mathcal{U})} \simeq
\bigwedge\limits^{q}\widehat{\mathcal O}_{11}\vert_{\mathcal{U}}\otimes \pi_{\mathcal S}^*  (\widehat{\mathcal
T}_{\mathcal{S}})_{p-q}$.
\end{enumerate} 
The tensor product $\otimes$ is taken over the corresponding sheaf of holomorphic functions. 

\end{proposition}

\medskip

\noindent {\it Proof.} The proof follows from the definitions.$\square$

\medskip

\section{Supermanifold of flags}\label{sec def Supermanifold of flags}

Denote by $\mathbf F^{m}_{k}$ the manifold of flags of type
$k=(k_0,k_1,\ldots,k_r)$ in $\mathbb C^m$, where $0 \leq k_r \leq\dots \leq
k_1 \leq k_0= m$. Let us describe an atlas on $\mathbf F^{m}_{k}$ that we will adapt for supercase. 
Let $\mathbb C^m\supset W_1\supset\dots\supset W_r$ be a flag of type
$k$. We choose a basis $B_s$ in each $W_s$ and assume that $B_0 = (e_1,\dots,e_m)$ is the standard basis in $\mathbb C^m$. Further, for any $s = 1,\ldots,r$, we define the matrix $X_s\in
\operatorname{Mat}_{k_{s-1},k_s}(\mathbb C)$  in the
following way: the columns of $X_s$ are the coordinates of the
vectors from $B_s$ with respect to the basis $B_{s-1}$. Since $\rk
X_s = k_s$, the matrix $X_s$ contains a non-degenerate minor of size
$k_s$.

For each $s = 1,\ldots,r$ let us fix a $k_s$-tuple
$I_s\subset\{1,\ldots,k_{s-1}\}$. We put $I = (I_1,\ldots,I_r)$. Denote
by $U_I$ the set of flags from $\mathbf F^{m}_{k}$ satisfying the
following conditions: there exist bases $B_s$ such that $X_s$
contains the identity matrix of size $k_s$ in the lines with numbers
from $I_s$. Clearly any flag from $U_I$ is uniquely
determined by those elements of $X_s$ that are not contained in the
identity matrix. Furthermore, any flag is contained in a certain
 $U_I$. The elements of $X_s$ that are not contained in the identity
 matrix are the coordinates of a flag
 from $U_I$ in the chart determined by $I$.
 Hence the local coordinates in $U_I$ are determined by
 $r$-tuple $(X_{1},\ldots,X_{r})$.

  Rename $X_{I_s}:= X_s$. If $J = (J_1,\ldots,J_s)$,
 where $J_s\subset\{1,\ldots,k_{s-1}\}$ and $|J_s| = k_s$, then the
 transition functions between the charts $U_I$ and $U_J$
 are given by:
$$
X_{J_1} = X_{I_1}C_{I_1J_1}^{-1}, \ \ X_{J_s} =
C_{I_{s-1}J_{s-1}}X_{I_s}C_{I_sJ_s}^{-1},\ \ s\ge 2,
$$
where $C_{I_1J_1}$ is the submatrix of $X_{I_1}$ formed by the lines
with numbers from $J_1$ and $C_{I_sJ_s}$, $s\ge 2$, is the submatrix
of $C_{I_{s-1}J_{s-1}}X_{I_s}$ formed by lines with numbers from
$J_s$.

Let us give a similar description of a classical flag supermanifold
in terms of atlases and local coordinates, see also\cite{ViGL in JA,
ViPi-sym}.  Let us take $m,n\in\mathbb N$ and
let 
$$
k=(k_0,k_1,\ldots,k_r)\quad \text{and} \quad l=(l_0,l_1,\ldots,l_r)
$$ 
be two $(r+1)$-tuples such that 
\begin{align*}
0\le k_r\le\ldots\le k_1\le k_0= m,\quad 0\le l_r\ldots\le
l_1\le l_0= n, \\  
0 < k_r +l_r <\ldots < k_1 + l_1 < k_0+l_0= m+n.
\end{align*}
 
 Let us define
the  supermanifold $\mathbf F^{m|n}_{k|l}$ of flags of type 
$(k|l)$ in the superspace  $V = \mathbb C^{m|n}$. The underlying manifold of
$\mathbf F^{m|n}_{k|l}$ is the product $\mathbf F^{m}_{k}
\times\mathbf F^{n}_{l}$ of two manifolds of flags of type $k$ and $l$
in $\mathbb C^m = V_{\bar 0}$ and $\mathbb C^n = V_{\bar 1}$.

For each $s = 1,\ldots,r$ let us fix $k_s$- and $l_s$-tuples of
numbers
 $I_{s\bar 0}\subset\{1,\ldots,k_{s-1}\}$ and $I_{s\bar
1}\subset\{1,\ldots,l_{s-1}\}$. We put
$I_s=(I_{s\bar 0},I_{s\bar 1})$ and $ I = (I_1,\ldots,I_r)$. Our goal now is to construct a superdomain $\mathcal{W}_I$. To each $I_s$ let us assign
a matrix of size $(k_{s-1} + l_{s-1})\times (k_s + l_s)$
\begin{equation}\label{Z_I_s}
\begin{split}
Z_{I_s}=\left(
\begin{array}{cc}
X_s & \Xi_s\\
\H_s & Y_s \end{array} \right), \ \ s=1,\dots,r.
\end{split}
\end{equation}
Assume that the identity matrix $E_{k_s+l_s}$ is contained in the
lines of $Z_{I_s}$ with numbers $i\in I_{s\bar 0}$ and $k_{s-1} +
j,\; j\in I_{s\bar 1}$. Here $X_s\in \operatorname{Mat}_{k_{s-1},
k_{s}}(\mathbb C),\; Y_s\in \operatorname{Mat}_{l_{s-1}, l_{s}}(\mathbb
C)$, where $\operatorname{Mat}_{a, b}(\mathbb C)$ is the space of
matrices of size $a\times b$ over $\mathbb C$. By definition, the
entries of $X_s$ and $Y_s$, $s=1,\dots,r$, that are not contained in
the identity matrix form the even coordinate system of
$\mathcal{W}_I$. The non-zero entries of $\Xi_s$ and $\H_s$ form the
odd coordinate system of $\mathcal{W}_I$.

Thus we have defined a set of superdomains on $\mathbf F^{m}_{k}
\times\mathbf F^{n}_{l}$ indexed by $I$. Note that the reductions of
these superdomains cover $\mathbf F^{m}_{k} \times\mathbf F^{n}_{l}$.
The local coordinates of each superdomain are determined by the
$r$-tuple of matrices $(Z_{I_1},\ldots,Z_{I_r})$. Let us define the
transition functions between two superdomains corresponding to $I =
(I_s)$ and $J = (J_s)$ by the following formulas:
\begin{equation}
\label{perehod}
\begin{split}
Z_{J_1} = Z_{I_1}C_{I_1J_1}^{-1}, \ \ Z_{J_s} =
C_{I_{s-1}J_{s-1}}Z_{I_s}C_{I_sJ_s}^{-1},\ \ s\ge 2,
\end{split}
\end{equation}
where $C_{I_1J_1}$ is the submatrix of $Z_{I_1}$ that consists of
the lines with numbers from $J_1$, and $C_{I_sJ_s}$, $s\ge 2$, is
the submatrix of $C_{I_{s-1}J_{s-1}}Z_{I_s}$ that consists of the
lines with numbers from $J_s$. Gluing the superdomains $\mathcal
W_I$, we define the {\it supermanifold of flags} $\mathbf
F^{m|n}_{k|l}$. In the case $r = 1$ this supermanifold is  called a
{\it super-grassmannian}. In the literature the notation $\mathbf
{Gr}_{m|n,k_1|l_1}$ is sometimes used.

The supermanifold $\mathbf F^{m|n}_{k|l}$ is $\GL_{m|n}(\mathbb
C)$-homogeneous. The action is given by
\begin{equation}
\label{dey}
\begin{split}
(L,(Z_{I_1},\ldots,Z_{I_r}))\mapsto (\tilde Z_{J_1},\ldots,\hat
Z_{J_r}),\\
\tilde Z_{J_1} = LZ_{I_1}C_1^{-1},\,\, \tilde Z_{J_s} =
C_{s-1}Z_{I_s}C_s^{-1}.
\end{split}
\end{equation}
Here $L$ is a coordinate matrix of $\GL_{m|n}(\mathbb C)$, $C_1$ is
the invertible submatrix of $LZ_{I_1}$ that consists of the lines
with numbers from $J_1$, à $C_s,\; s\ge 2$, is the invertible
submatrix of $C_{s-1}Z_{I_s}$ that consists of the lines with
numbers from $J_s$. If $0<k_r<\ldots <k_1<m,\quad 0<l_r<\ldots <l_1<n$ we will say that the flag supermanifold $\mathbf F^{m|n}_{k|l}$ has {\it generic type}.

\section{Retract of the superbundle $\mathbf F^{m|n}_{k|l}$}\label{sec retract of flag smf}

Let $\mathcal G$ be a Lie supergroup, $\mathcal H^1$ and $\mathcal H^2$ be two
 Lie subsupergroups of $\mathcal G$ such that $\mathcal H^2$ is a Lie subsupergroup of
 $\mathcal H^1$. (More information about Lie supergroups can be found for instance in \cite{V_funk,ViLieSupergroup}.) Then the homogeneous superspace $\mathcal M = \mathcal G/\mathcal H^2$ is a
 superbundle with base $\mathcal B = \mathcal G/\mathcal H^1$ and fiber $\mathcal F = \mathcal H^1/\mathcal H^2$.
In \cite{V_funk} it was proved that 
$$
\gr\mathcal M = \gr \mathcal G/ \gr \mathcal H^2, \quad
\gr\mathcal B = \gr \mathcal G/ \gr \mathcal H^1.
$$ 
Moreover, $\gr \mathcal H^2 $ is a Lie subsupergroup of $\gr \mathcal H^1$. The superbundle
$\tilde{\pi}:\gr \mathcal G/ \gr \mathcal H^2 \to \gr \mathcal G/ \gr \mathcal H^1$,
which we will also denote by $\tilde{\pi}:\tilde{\mathcal M}\to
\tilde{\mathcal B}$,  is the split superbundle corresponding to $\mathcal M\to \mathcal B$, see Subsection \ref{sec Retract of a split superbundle}.
Denote by
 ${\pi}: \widehat {\mathcal M}\to
\widehat {\mathcal B}$ the retract of
$\tilde{\pi}:\tilde{\mathcal M}\to
\tilde{\mathcal B}$. As above we denote by $\mathcal H^1_0$ and $\mathcal H^2_0$ the underlying spaces of $\mathcal H^1$ and $\mathcal H^2$, respectively.

Let $G$ be a Lie group and $H$ be its Lie subgroup. Consider a locally free sheaf $\mathcal E$ on the homogeneous manifold $G/H$ and denote by $\mathbb E$ the corresponding vector bundle. (Recall that there is a one-to-one correspondence between locally free sheaves on $G/H$ and vector bundles over $G/H$.) Assume that $\mathbb E$ is a homogeneous vector bundle, see \cite{ADima}. In this case sometimes we will call $\mathcal E$ a homogeneous locally free sheaf. Recall that there is a bijection between homogeneous bundles $\mathbb E$ over $G/H$ and representations of $H$ in the fiber $\mathbb E_{eH}$, see \cite{ADima}. For simplicity sometimes we will write: the representation of $H$ corresponding to the homogeneous locally free sheaf $\mathcal E$ meaning this correspondence.

Denote by $\alpha$ the representation of the Lie group $\mathcal H^1_0$ corresponding to the homogeneous locally free sheaf $(\tilde{\mathcal{O}}_{\mathcal B})_1$ and by  $\beta$ the representation of $\mathcal H^2_0$ corresponding to the homogeneous locally free sheaf $\tilde{\mathcal{O}}_1$. (As above we denote by $\tilde{\mathcal{O}}$ and by $\widehat{\mathcal{O}}$ the structure sheaf of $\tilde{\mathcal{M}}$ and $\widehat{\mathcal{M}}$, respectively.) Recall that
$$
\widehat{\mathcal{O}}_1 = \widehat{\mathcal{O}}_{11} + \widehat{\mathcal{O}}_{10} = {\pi}^*((\mathcal{O}_{\mathcal B})_1) + (\tilde{\mathcal{O}})_1/{\pi}^*((\mathcal{O}_{\mathcal B})_1).
$$
Hence the locally free sheaf $\widehat{\mathcal{O}}_1$ corresponds
to the representation 
\begin{equation}\label{eq represent in tilda O_1 alpha}
\alpha|_{\mathcal H^2_0} + \beta/ (\alpha|_{\mathcal H^2_0})
\end{equation}
of Lie group $\mathcal H^2_0$. In particular,  $\tilde{\mathcal M}$ is a
homogeneous supermanifold. Note that a split supermanifold is homogeneous if and only if the corresponding bundle is homogeneous, see \cite{V_funk}.

Let $\mathcal M: = \mathbf F^{m|n}_{k|l}$ and $r\geq 2$. The flag supermanifold $\mathcal M$ is a superbundle with base
$\mathcal B: = \mathbf F^{m|n}_{k_1|l_1}$ and fiber $\mathcal S: = \mathbf F^{k|l}_{k'|l'}$, where
\begin{equation}\label{eq def of k'= l'=}
k'= (k_1,\ldots,k_r)\quad \text{and} \quad l'=(l_1,\ldots,l_r).
\end{equation} 
In coordinates (\ref{Z_I_s}) the bundle projection $\pi: \mathcal M \to \mathcal B$ is given by 
$$
(Z_{I_1},\ldots,Z_{I_r}) \longmapsto (Z_{I_1}).
$$

Denote by $\tilde{\pi}:\tilde{\mathcal M}\to \tilde{\mathcal B}$
the corresponding to $\mathbf F^{m|n}_{k|l}$ split superbundle with fiber
$\tilde{\mathcal S}$, see Subsection \ref{sec Retract of a split superbundle}. Denote also by
${\pi}:\widehat {\mathcal M}\to
\widehat {\mathcal B}$ the retract of
$\tilde{\pi}:\tilde{\mathcal M}\to \tilde{\mathcal B}$.
 Consider the superdomain $Z_{I}$ in $\mathbf F^{m|n}_{k|l}$ corresponding to
 \begin{equation}\label{eq I_0 and I_1 of Z_I}
 I_{s\bar 0} = (k_{s-1}-k_s+1,\ldots,k_{s-1}),\,\,\,I_{s\bar 1} =
 (l_{s-1}-l_s+1,\ldots,l_{s-1}),
 \end{equation}
 where $s=1,\ldots, r$, see Section \ref{sec def Supermanifold of flags}. Denote by $x$ the
origin of $Z_{I}$. We also denote by $P$ and by $\mathfrak p$ the underlying Lie group and Lie superalgebra of the
super-stabilizer of $x$ for the action (\ref{dey}) of
$\GL_{m|n}(\mathbb{C})$. Let $R$ be the reductive part of $P$. A direct calculation shows that $R$ has the following form
\begin{equation}\label{eq reductive part R=}
R =\left\{ \left(
    \begin{array}{ccc}
      A_1 & 0 & 0 \\
      0& A_2 & \vdots \\
      0 & \cdots & \ddots \\
    \end{array}
  \right)\times
\left(
    \begin{array}{ccc}
      B_1 & 0 & 0 \\
      0 & B_2 & \vdots \\
      0& \cdots & \ddots \\
    \end{array}
  \right) \right\},
\end{equation}
where $A_i\in \GL_{k_{i-1}-k_i} (\mathbb{C})$ and $B_i\in
\GL_{l_{i-1}-l_i} (\mathbb{C})$, where $k_{r+1}=l_{r+1}:=0$ and $i=1,\ldots,
r+1$. Denote by $\rho_i$ the standard representations of the group
$\GL_{k_{i-1}-k_i} (\mathbb{C})$ and by $\sigma_i$ the standard
representations of the Lie group $\GL_{l_{i-1}-l_i} (\mathbb{C})$.

The sheaf $\tilde{\mathcal{O}}_1$ is a homogeneous locally free
sheaf on $\GL_{m}(\mathbb{C})\times \GL_{n}(\mathbb{C})/P$. Let us compute the corresponding representation $\theta$ of
$P$. The fiber over $x$ of the
corresponding to $\tilde{\mathcal{O}}_1$ vector bundle is isomorphic to  $(\mathfrak{gl}_{m|n}(\mathbb C)_{\bar 1}/\mathfrak{p}_{\bar 1} )^*$ as
$P$-modules, see \cite{V_funk} for details. A direct computation shows
that
$$
\mathfrak{p}_{\bar 1} = 
\left\{ \left(
                                               \begin{array}{cc}
                                                     0 & C \\
                                                     D & 0 \\
                                                   \end{array}\right)\,\,
                                                   |\,\,\, C\in
                                                   \Mat_{m,n}(\mathbb{C}),
\,\,\, D\in
                                                   \Mat_{n,m}(\mathbb{C})
 \right\} \subset
\mathfrak{gl}_{m|n}(\mathbb C)_{\bar 1},
$$
where
$$
C=
\left(
    \begin{array}{ccc}
      C_1 & 0 & 0 \\
      *& C_2 & \vdots \\
      * & \cdots & \ddots \\
    \end{array}
  \right),  \,\,\,\,\,
  D=
\left(
    \begin{array}{ccc}
      D_1 & 0 & 0 \\
      *& D_2 & \vdots \\
      * & \cdots & \ddots \\
    \end{array}
  \right),
$$
where $C_i\in \Mat_{k_{i-1}- k_{i}, l_{i-1}- l_{i}}(\mathbb{C})$ and
$D_i\in \Mat_{l_{i-1}- l_{i}, k_{i-1}- k_{i}}(\mathbb{C})$.

Denote by $\varphi$ the representation of $P$
corresponding to the locally free sheaf $\widehat{\mathcal{O}}_{11}$
and by $\psi$ the representation of $P$
corresponding to the locally free sheaf
$\widehat{\mathcal{O}}_{10}$.

\begin{lemma}\label{lemma_reps of O_11 nd O_10} 
We have
$$
\varphi|R= \sum_{i>1} \rho_1^*\otimes\sigma_{i}+
\sum_{i>1}\sigma_1^*\otimes\rho_{i},
\quad
\psi|R= \sum_{1<i<j} \rho_i^*\otimes\sigma_{j}+
\sum_{1<i<j}\sigma_i^*\otimes\rho_{j},
$$
$$
\theta|R = \varphi|R+
\psi|R = \sum_{i<j} \rho_i^*\otimes\sigma_{j}+
\sum_{i<j}\sigma_i^*\otimes\rho_{j}.
$$
\end{lemma}

\medskip

\noindent {\it Proof.} Let us compute the representation
$\theta^*|R$ of $R$ in
$\mathfrak{gl}_{m|n}(\mathbb C)_{\bar 1}/
\mathfrak{p}_{\bar 1}$ for example for $r=1$.
We identify the vector space $\mathfrak{gl}_{m|n}(\mathbb C)_{\bar 1}/
\mathfrak{p}_{\bar 1}$ with 
$$
\mathfrak{p}_{\bar 1} = 
\left\{ \left(
\begin{array}{cc}
0 &\tilde C \\
\tilde D & 0 \\
\end{array}\right)\,\,
|\,\,\, \tilde C\in
\Mat_{m,n}(\mathbb{C}),
\,\,\, \tilde D\in
\Mat_{n,m}(\mathbb{C})
\right\} \subset
\mathfrak{gl}_{m|n}(\mathbb C)_{\bar 1},
$$
where
$$
\tilde C=\left( \begin{array}{cc}
0 & M \\
0 & 0 \\
\end{array}\right) \quad \text{and} \quad \tilde D=\left( \begin{array}{cc}
0 & N \\
0 & 0 \\
\end{array}\right),
$$
 $M\in \Mat_{k_0-k_1,l_1-l_2}(\mathbb{C})$ and $N\in \Mat_{l_0-l_1,k_1-k_2}(\mathbb{C})$. We have

$$
\begin{array}{c}
\left[\left(
\begin{array}{cccc}
A_1 & 0 & 0 & 0 \\
0 & A_2 & 0 & 0 \\
0 & 0 & B_1 & 0 \\
0 & 0 & 0 & B_2 \\
\end{array}
\right), \left(
\begin{array}{cccc}
0 & 0 & 0 & M \\
0 & 0 & 0 & 0 \\
0 & N & 0 & 0 \\
0 & 0 & 0 & 0 \\
\end{array}
\right)
\right] = \\
\left(
\begin{array}{cccc}
0 & 0 & 0 & A_1M - MB_2 \\
0 & 0 & 0 & 0 \\
0 & B_1N - NA_2 & 0 & 0 \\
0 & 0 & 0 & 0 \\
\end{array}
\right).
\end{array}
$$
Hence, $\theta^*|R = \rho_1\otimes \sigma_2^* +
\sigma_1\otimes \rho_2^*$. For $r\geq 2$ the proof is
similar. To obtain $\varphi|R$ and $\psi_R$ we use Formula (\ref{eq represent in tilda O_1 alpha}).$\Box$

\medskip

Recall that we denoted by $\Theta$ and by $\Theta^v$ the sheaf of vector fields and the sheaf of vertical vector fields on the manifold $\mathcal M_0$, respectively. 
Denote by $\tau$, $\tau^v$ and $\tau^h$ the
representation of $P$ corresponding to the sheaves $\Theta$,  $\Theta^v$ and $\Theta^h=\Theta/\Theta^v$, respectively. (Recall that we have a natural representation of $P$ in the fiber over the
origin $x$ of $Z_{I}$.)

\begin{lemma} \label{lemma_reps of T,T^v and T/T^v}   We have
$$
\begin{array}{c}
\tau|R = \sum\limits_{i<j}
(\rho_i\otimes\rho_j^*+\sigma_i\otimes\sigma_j^*),\,\,\,
\tau^v|R = \sum\limits_{1<i<j}
(\rho_i\otimes\rho_j^*+\sigma_i\otimes\sigma_j^*),\\
\tau^h|R = \sum\limits_{1<j}
(\rho_1\otimes\rho_j^*+\sigma_1\otimes\sigma_j^*).
\end{array}
$$
\end{lemma}

\medskip

\noindent {\it Proof.}
The representation $\tau|R$ is the isotropy representation. Recall that the isotropy representation is isomorphic to the natural representation of $P$ in the vector space $\mathfrak{gl}_{m|n}(\mathbb C)_{\bar 0}/\mathfrak p_{\bar 0}$. Let us compute this representation in the case $r=2$. We identify $\mathfrak{gl}_{m|n}(\mathbb C)_{\bar 0}/\mathfrak p_{\bar 0}$  with
$$
\left\{ 
\left(
\begin{array}{ccc}
* & U_{12} & U_{13}  \\
* & * & U_{23}  \\
* & * & *  \\
\end{array}
\right) \times
\left(
\begin{array}{ccc}
* & V_{12} & V_{13}  \\
* & * & V_{23}  \\
* & * & *  \\
\end{array}
\right)
\right\},
$$
where $U_{ij} \in \Mat_{k_{i-1}-k_i,k_{j-1}-k_j}(\mathbb{C})$ and $V_{ij} \in \Mat_{l_{i-1}-l_i,l_{j-1}-l_j}(\mathbb{C})$. 
 We have
$$
\left[
\left(
\begin{array}{ccc}
A_1 & 0 & 0  \\
0 & A_2 & 0  \\
0 & 0 & A_3  \\
\end{array}
\right),
\left(
\begin{array}{ccc}
* & U_{12} & U_{13}  \\
* & * & U_{23}  \\
* & * & *  \\
\end{array}
\right)
\right]= \left(
\begin{array}{ccc}
* & A_1U_{12}-U_{12}A_2 & A_1U_{13}-U_{13}A_3  \\
* & * & A_2U_{23}-U_{23}A_3  \\
* & * & *  \\
\end{array}
\right),
$$
$$
\left[
\left(
\begin{array}{ccc}
B_1 & 0 & 0  \\
0 & B_2 & 0  \\
0 & 0 & B_3  \\
\end{array}
\right),
\left(
\begin{array}{ccc}
* & V_{12} & V_{13}  \\
* & * & V_{23}  \\
* & * & *  \\
\end{array}
\right)
\right]= \left(
\begin{array}{ccc}
* & A_1V_{12}-V_{12}A_2 & A_1V_{13}-V_{13}A_3  \\
* & * & A_2V_{23}-V_{23}A_3  \\
* & * & *  \\
\end{array}
\right).
$$
Therefore 
$$
\tau|R= \sum_{i<j} (\rho_i\otimes\rho_j^*+\sigma_i\otimes\sigma_j^*),
$$
For $\tau^v|R$ and $\tau^h|R$ we respectively have
$$
\tau^v|R= \rho_2\otimes\rho_3^*+\sigma_2\otimes\sigma_3^*,\,\,\,\tau^h|R= \rho_1\otimes\rho_2^* +\rho_1\otimes\rho_3^*+
\sigma_1\otimes\sigma_2^*+\sigma_1\otimes\sigma_3^*.
$$
 For $r\geq 2$ the proof is similar.
$\Box$

\medskip

\section{ Vector fields on the retract of flag super\-manifold}

Recall that $\mathcal M = \mathbf F^{m|n}_{k|l}$, where $r\geq 2$, is a
superbundle with base space $\mathcal B = \mathbf F^{m|n}_{k_1|l_1}$ and fiber $\mathcal S = \mathbf F^{k_1|l_1}_{k'|l'}$, where $k',l'$ are defined in (\ref{eq def of k'= l'=}). As above we denote by
$\tilde{\pi}:\tilde{\mathcal M}\to \tilde{\mathcal B}$
 the corresponding to $\mathcal M$ split superbundle with fiber
$\gr\mathcal S$, see Section \ref{sec retract of flag smf}. The aim of this
section is to compute the Lie superalgebra of holomorphic vector
fields on the retract $\tilde{\mathcal M}$ of the flag
supermanifold $\mathcal M$ for $r\geq 2$.

In Section \ref{sec Split supermanifolds} we defined the functor $\gr$. Recall that $\mathbf{Gr}_{m\vert n, k\vert l}$ is the super-grass\-mannian of type $(k|l)$. 
We denote the Lie superalgebra of holomorphic vector fields on $\gr\mathbf{Gr}_{m\vert n, k\vert l}$ by $\mathfrak
v(\gr\mathbf{Gr}_{m\vert n, k\vert l})$. Consider the following Cartan
subalgebra in the Lie algebra $\mathfrak {gl}_{m\vert n}(\mathbb{C})_{\bar 0}$
$$
\mathfrak{h} := \{\operatorname{diag}(\mu_1,\dots,\mu_{m})\}\oplus
\{\operatorname{diag}(\lambda_1,\dots,\lambda_{n})\}.
$$ 
The Lie superalgebra  $\mathfrak
v(\gr\mathbf{Gr}_{m\vert n, k\vert l})$ was computed in \cite{onigl}.

\begin{theorem}\cite[Theorem 4, Lemma 3]{onigl}\label{theor vect fields on Gr}
\label{teor vector fields Grassmannian gr}  Let $0<k<
	m$ and $0<l< n$. Then 
	$$
	\mathfrak v(\gr\mathbf{Gr}_{m\vert n, k\vert
		l})= \mathfrak v(\gr\mathbf{Gr}_{m\vert n, k\vert l})_{-1} \oplus
	\mathfrak v(\gr\mathbf{Gr}_{m\vert n, k\vert l})_{0},
	$$ 
	where 
	$$
	\mathfrak v(\gr\mathbf{Gr}_{m\vert n, k\vert l})_{-1} \simeq \mathfrak {gl}_{m\vert n}(\mathbb C)_
	{\bar 1}
	$$ 
	as $\mathfrak {gl}_{m\vert n}(\mathbb{C})_{\bar
		0}$-modules. For $\mathfrak v(\gr\mathbf{Gr}_{m\vert n, k\vert l})_{0}$  we have the following possibilities. 
	\begin{enumerate}
		
		\item  $\mathfrak v(\gr\mathbf{Gr}_{m\vert n, k\vert l})_{0}\simeq \mathfrak {gl}_{m\vert n}(\mathbb{C})_{\bar
			0}$ (as Lie algebras) for 
		{\bf (1)} $1<k<m-1$ and $1< l< n-1$; {\bf (2)} $k=1$ and $l< n-1$; {\bf (3)} $k=m-1$ and $l>1$; {\bf (4)} $k<m-1$ and $l=1$; {\bf (5)} $k>1$ and $l= n-1$.

		\item if {\bf (1)} $k=n-l=1$ and $m-k>1$ or {\bf (2)} $k=n-l=1$ and  $l>1$, then 
		$$
		\mathfrak
		v(\gr\mathbf{Gr}_{m\vert n, k\vert l})_{0}\simeq \mathfrak
		{gl}_{m\vert n}(\mathbb{C})_{\bar 0} \oplus \mathfrak k
		$$ (as
		$\mathfrak {gl}_{m\vert n}(\mathbb{C})_{\bar 0}$-modules), where
		$\mathfrak k$ is the irreducible $\mathfrak {gl}_{m\vert
			n}(\mathbb{C})_{\bar 0}$-module with the highest weight
		$-\mu_{m-1}-\mu_{m}+\lambda_{1}+\lambda_{2}$;
		
		\item if {\bf (1)} $l=m-k=1$ and $n-l>1$, or {\bf (2)} $l=m-k=1$ and $k>1$, then 
		$$
		\mathfrak
		v(\gr\mathbf{Gr}_{m\vert n, k\vert l})_{0}\simeq \mathfrak
		{gl}_{m\vert n}(\mathbb{C})_{\bar 0}\oplus \mathfrak k
		$$ 
		(as
		$\mathfrak {gl}_{m\vert n}(\mathbb{C})_{\bar 0}$-modules), where
		$\mathfrak k$ is the irreducible $\mathfrak {gl}_{m\vert
			n}(\mathbb{C})_{\bar 0}$-module with the highest weight
		$\mu_{1}+\mu_{2}-\lambda_{n-1}-\lambda_{n}$;

		\item if $m=n=2$, $k=l=1$, then 
		$$
		\mathfrak v(\gr\mathbf{Gr}_{m\vert n, k\vert l})_{0}\simeq \mathfrak {gl}_{m\vert n}(\mathbb{C})_{\bar
			0}\oplus \mathfrak k\oplus \mathfrak k'
		$$ 
		(as $\mathfrak {gl}_{m\vert
			n}(\mathbb{C})_{\bar 0}$-modules), where $\mathfrak k$, $\mathfrak
		r'$ are the irreducible $\mathfrak {gl}_{m\vert n}(\mathbb{C})_{\bar
			0}$-modules with the highest weights
		$-\mu_{1}-\mu_{2}+\lambda_{1}+\lambda_{2}$ and
		$\mu_{1}+\mu_{2}-\lambda_{1}-\lambda_{2}$ respectively;
		
	\end{enumerate}
\end{theorem}

Note that in \cite[Theorem 4]{onigl} there is a mistake in the statement. To obtain the result of Theorem \ref{teor vector fields Grassmannian gr} we need to use \cite[Theorem 4]{onigl} and \cite[Lemma 3]{onigl}. Further, the case $k=0$, when the super-grassmannian is split, was considered by A.~Onishchik and A.~Serov. A complete description of the Lie superalgebras of holomorphic vector fields on super-grassmannians can be found in \cite{ViGL in JA}. The Lie superalgebra of holomorphic vector fields on a flag
supermanifold $\mathcal M=\mathbf F^{m|n}_{k|l}$ of length $r>1$ was computed in generic case in \cite{ViGL in JA}. Using the method developed in \cite{ViGL in JA} we can calculate the Lie superalgebra of holomorphic vector fields on $\tilde{\mathcal M}= \gr \mathbf F^{m|n}_{k|l}$. Since this computation is similar we will omit some details.

In what follows we will consider flag supermanifolds under the following
conditions on the flag type:
\begin{equation}\label{eq generic type condition for flags}
0< k_r < \cdots < k_1< k_0=m\quad \text{and} \quad 0< l_r < \cdots < l_1< l_0=n
\end{equation}
and
\begin{equation}\label{condition flag type 1}
\begin{array}{rl}
(k_r,l_r)\ne &(1,l_{r-1}-1), (k_{r-1}-1,1), (1,l_{r-1}-2),\\
&(k_{r-1}-2,1), (2,l_{r-1}-1), (k_{r-1}-1,2).
\end{array}
 \end{equation}
Assumption (\ref{eq generic type condition for flags}) is related to the fact that under this conditions $\mathfrak v({\mathcal M})\simeq \mathfrak {gl}_{m\vert n}(\mathbb C)$, see \cite{ViGL in JA}. We assume (\ref{condition flag type 1}) due to our use of induction and results of \cite{onigl}. Indeed, results \cite{onigl} were obtained under conditions (\ref{condition flag type 1}) for a super-grassmannian.

Denote by $\mathfrak v(\tilde{\mathcal M})$ the Lie superalgebra of holomorphic vector fields on $\tilde{\mathcal M}$. 

\begin{theorem}\label{theor vect fields on flag}
	 Assume that $r>1$ and $\mathbf F^{m|n}_{k|l}$ satisfies the conditions (\ref{eq generic type condition for flags}) and (\ref{condition flag type
 1}). Then $\mathfrak v(\tilde{\mathcal M})=
\mathfrak v(\tilde{\mathcal M})_{-1}\oplus \mathfrak
v(\tilde{\mathcal M})_{0}$, where
\begin{enumerate}
\item $\mathfrak v(\tilde{\mathcal M})_{-1} \simeq \mathfrak {gl}_{m\vert n}(\mathbb C)_{\bar 1}$ (as
$\mathfrak {gl}_{m\vert n}(\mathbb C)_{\bar 0}$-modules);
\item $\mathfrak v(\tilde{\mathcal M})_{0}\simeq \mathfrak {gl}_{m\vert n}(\mathbb C)_{\bar 0}$ (as Lie algebras).
\end{enumerate}
\end{theorem}

\noindent {\it Proof.}   In \cite{V_funk} it was shown that
$H^0(\mathcal M_0,\tilde{\mathcal O})\simeq \mathbb C$ if (\ref{eq generic type condition for flags}) hold true.   Hence by Theorem \ref{bash} any vector field on
$\tilde{\mathcal M}$ is projectible with respect to
$\tilde{\pi}$. It follows that we have a homomorphism $\Pi:
\mathfrak v(\tilde{\mathcal M})\to \mathfrak v(\tilde{\mathcal B})$ of Lie superalgebras. The idea of \cite{ViGL in JA} was to compute the kernel and the
image of $\Pi$.  

Denote by $\mathcal W =
(\tilde\pi_{0})_*\mathcal T^v$ the direct image of the sheaf of
vertical vector fields on $\tilde{\mathcal M}$. It is clear that $\Ker\Pi =
H^0(\mathcal M_0,\mathcal{W})$.  The sheaf $\mathcal W$ possesses the following
filtration:
$$
\mathcal W= \mathcal W_{(0)} \supset \mathcal W_{(1)} \supset
\ldots,
$$
where $\mathcal W_{(p)}:= (J_{\mathcal B})^p \mathcal{W}$ and $J_{\mathcal B}$ is the ideal in $\tilde{\mathcal O}$ generated by odd elements in $\tilde\pi^*(\tilde{\mathcal O}_{\mathcal B})$. We put
$\widetilde{\mathcal W}_p:= \mathcal W_{(p)}/\mathcal W_{(p+1)}$. The
sheaves $\widetilde{\mathcal W}_p$ are homogeneous locally free sheaves of $\mathcal{F}_{\mathcal B_0}$-modules, see \cite{ViGL in JA} for details. The $0$-cohomology
group of $\widetilde{\mathcal W}_p$ can be computed using the Borel -
Weil - Bott Theorem. Denote by $P_{\mathcal B}$ the stabilizer of the origin of the chart $Z_{I_1}$, see (\ref{eq I_0 and I_1 of Z_I}) and below.  A direct computation shows that the dominant highest weight of the representation of $P_{\mathcal B}$ corresponding to $\mathcal W_0$ is $0\,\,(\times 2)$. Therefore by the Borel - Weil - Bott Theorem $H^0(\mathcal B_0, \widetilde{\mathcal W}_0)= \mathbb C^2$. Similarly to \cite{ViGL in JA} we can prove that $H^0(\mathcal B_0, \mathcal W_{(0)})= 0$. The idea here is to use \cite[Lemma 2]{ViGL in JA} and the fact that $H^0(\mathcal B_0, \mathcal W_{(0)})$ is an ideal in the Lie superalgebra of vector fields on  $\tilde{\mathcal M}$. Therefore the homomorphism $\Pi:
\mathfrak v(\tilde{\mathcal M})\to \mathfrak v(\tilde{\mathcal B})$ is injective.

Let us show that $\Pi$ is surjective. 
By Theorem  \ref{theor vect fields on Gr} on the base space we have  $\mathfrak v(\tilde {\mathcal B}) \simeq \mathfrak {gl}_{m\vert n}(\mathbb C)$ (as $\mathfrak {gl}_{m\vert n}(\mathbb C)_{\bar 0}$-modules). Further we have the following embedding $\mathfrak {pgl}_{m\vert n}(\mathbb C)=\mathfrak v(\mathcal M) \hookrightarrow \mathfrak v(\tilde{\mathcal M})$, which image is equal to the Lie superalgebra of all fundamental vector fields on $\tilde{\mathcal M}$. Moreover on any split supermanifold there is the grading operator which we denote by $z$. This operator acts in the sheaf $\tilde{\mathcal O}$ in the following way: $ z(f)= pf$, where $f\in \tilde{\mathcal O}_p$. The operator $z$ is not fundamental since there is no such operators on  $\mathcal M$. Further by definition of $z$ we have $[z, \tilde{\mathcal T}_o]=0$. Hence, 
$$
\mathfrak v(\tilde{\mathcal M})_{0}\simeq \mathfrak {pgl}_{m\vert n}(\mathbb C)_{\bar 0}\oplus \langle z\rangle \simeq \mathfrak {gl}_{m\vert n}(\mathbb C)_{\bar 0}. 
$$ 

By Theorem \ref{teor vector fields Grassmannian gr} under the conditions of our theorem all vector fields on $\tilde{\mathcal B}$ are fundamental with respect to the action of $\gr (\GL_{m|n})$ exept of the grading operator.  Therefore all holomorphic vector fields on $\tilde{\mathcal B}$ can be lifted to $\tilde{\mathcal M}$. In other words the homomorphism $\Pi$ is surjective.$\square$

\medskip
As a corollary we get.

\begin{corollary}\label{cor H^0(A_p)=}
Assume that $r>0$ and (\ref{eq generic type condition for flags}) and (\ref{condition flag type
	1}) hold true. Then
$$
H^0(\mathcal M_0,\mathcal A_p)= \left\lbrace
\begin{array}{l}
\lbrace0\rbrace, \,\,\,p\ne 0,-1;\\
\mathbb C^2, \,\,\,p= 0;\\
\mathfrak {gl}_{m\vert n}(\mathbb C)_{\bar 1}, \,\,\,p= -1,
\end{array}\right.
$$
where $\mathcal A_p = \tilde{\mathcal O}_1^*\otimes \bigwedge^p
\tilde{\mathcal O}_1$.

\end{corollary}

\noindent {\it Proof.} Consider exact Sequence (\ref{exact
sequens for T_p}). It determines the exact sequence of cohomology
groups
$$
 0\to H^0(\mathcal M_0,\mathcal A_{p}) \to H^0(\mathcal M_0,\tilde{\mathcal T}_{p})\to
H^0(\mathcal M_0,\mathcal C_{p}), \,\,\,p\geq -1.
$$
Using Theorem \ref{theor vect fields on flag} we conclude that
$$
H^0(\mathcal M_0,\mathcal A_{p})=\{0\}, \quad p\ne 0,-1. 
$$ 
 Further for $p= -1$ we have $\tilde{\mathcal T}_{-1}\simeq \mathcal A_{-1}$, hence
$$
H^0(\mathcal M_0,\mathcal A_{-1})\simeq H^0(\mathcal M_0,\tilde{\mathcal T}_{-1})\simeq \mathfrak {gl}_{m\vert n}(\mathbb C)_{\bar 1}.
$$ 
Consider the
case $p= 0$. It is a classical result  that the Lie algebra of holomorphic vector fields $H^0(\mathcal M_0,\Theta)$ on $\mathcal M_0$ is isomorphic to $\mathfrak
{sl}_{m}(\mathbb{C}) \oplus \mathfrak {sl}_{n}(\mathbb{C})$, see for example \cite{ADima} for details. Hence, 
$$
H^0(\mathcal M_0,\mathcal C_{0})= H^0(\mathcal M_0,\Theta)\simeq
\mathfrak {sl}_{m}(\mathbb{C})\oplus \mathfrak
{sl}_{n}(\mathbb{C}).
$$ 
The following diagram is commutative
$$
\begin{array}{ccc}
H^0(\mathcal M_0,\tilde{\mathcal T}_{0})& \longrightarrow  &H^0(\mathcal M_0,\Theta)\\
\downarrow && \downarrow\\
H^0(\mathcal B_0,(\tilde{\mathcal T}_{\mathcal B_0})_{0})& \longrightarrow &H^0(\mathcal B_0,\Theta_{\mathcal B_0})
\end{array}.
$$
The vertical right arrow is an isomorphism of Lie algebras since these Lie algebras are both isomorphic to $\mathfrak {sl}_{m}(\mathbb{C}) \oplus \mathfrak {sl}_{n}(\mathbb{C})$ and they coincide with Lie algebras of fundamental vector fields. The vertical left arrow is an isomorphism of Lie algebras by Theorem \ref{theor vect fields on Gr} and Theorem \ref{theor vect fields on flag}. Further, the map $H^0(\mathcal B_0,(\tilde{\mathcal T}_{\mathcal B_0})_{0}) \rightarrow H^0(\mathcal B_0,\Theta_{\mathcal B_0})$ is surjective, see \cite{onigl}.  It follows that the map $H^0(\mathcal M_0,\tilde{\mathcal T}_{0})\to
H^0(\mathcal M_0,\mathcal C_{0})$ is surjective too. Our assertion follows from
Theorem \ref{theor vect fields on flag}. $\square$
\smallskip

\section{Some known results about vector bundles}

Let us recall some results about vector bundles over usual complex-analytic manifolds. Let $X$ and $Y$ be complex-analytic manifolds, $\mathcal G$ and $\mathcal H$ be coherent analytic sheaves on $X$ and $Y$, respectively. Assume in addition that $X$ is a Stein manifold and $Y$ is compact. By the Cartan Theorems A and B we have $H^p(X,\mathcal G)=0$ for $p>0$ and by the Cartan-Serre Theorem $H^p(Y,\mathcal H)$ is finite dimensional. Now consider the direct product $X\times Y$ and the natural projections $\pi_X: X\times Y \to X$ and $\pi_Y: X\times Y \to Y$. As usual we denote by $\pi_X^*(\mathcal G)$ and by $\pi_Y^*(\mathcal H)$  the  pullback sheaves. These sheaves are sheaves of $\mathcal F_{X\times Y}$-modules, where $\mathcal F_{X\times Y}$ is the sheaf of holomorphic functions on $X\times Y$. The following theorem was proved in \cite{Kaup}. 

\begin{theorem}[K\"{u}nneth formula]\label{theor Kunneth Formula}
	Assume that $X$ is Stein, $Y$ is compact, $\mathcal G$ and $\mathcal H$ are locally free sheaves on $X$ and $Y$, respectively. 
	Then we have the following isomorphism
	\begin{equation}\label{eq Kunneth Formula}
H^k(X\times Y, \pi_X^*(\mathcal G)\otimes_{\mathcal F_{X\times Y}} \pi_Y^*(\mathcal H))	\simeq 	 H^0(X,\mathcal G) \otimes H^k(Y,\mathcal H). \Box
	\end{equation}
\end{theorem}

Let $\pi: M\to B$ be a complex-analytic bundle with compact fiber $F$ and $\mathbb{E}\to B$ be a holomorphic vector bundle. We denote by $\Gamma (V, \mathbb{E})$ the vector space of all sections of $\mathbb E$ over an open set $V\subset B$ and by $\pi^*\mathbb{E}\to M$ the pullback bundle.

\begin{lemma}\label{lem sections with compact fiber}   For any open set $V\subset B$ we have
	$$
	\Gamma(\pi^{-1}(V),\pi^*\mathbb{E}\vert_{\pi^{-1}(V)})\simeq \Gamma(V,\mathbb{E}\vert_{V} ).
	$$
\end{lemma}

\noindent{\it Proof.} 
The result is a consequence of the fact that 
\begin{equation}\label{eq pul back of vector bundles}
H^0(\pi^{-1}(V),\pi^*(\mathcal F_B))= H^0(V,\mathcal F_B),
\end{equation}
where $\mathcal F_B$ is the sheaf of holomorphic functions over $B$. Formula (\ref{eq pul back of vector bundles}) follows from Theorem \ref{theor Kunneth Formula}.$\Box$

\medskip

We will need also the following lemma that we prove here for completeness. 

\begin{lemma}\label{lem H^1(covering) = H^1}
	Let $\mathcal G$ be a sheaf of $\mathcal F_X$-modules on a complex-analytic manifold $X$ and $\mathcal U= \lbrace U_i\rbrace$ be an open covering of $X$. Assume that $H^1(U_i,\mathcal G)=\lbrace0\rbrace$ for any $U_i\in \mathcal U$. Then 
	$$
	H^1(X,\mathcal G )=H^1(\mathcal U, \mathcal G).
	$$
	In other words in this case we can compute 1st cohomology of $\mathcal G$ using the open covering $\mathcal U$. 
	
\end{lemma}

\noindent {\it Proof.} We put
$$
\mathcal H^p:= \prod_{i_0,\ldots,i_p}j_*(\mathcal G\vert_{ U_{i_0,\ldots,i_p}}), \quad p\geq 0,
$$
where $U_{i_0,\ldots,i_p}:=  U_{i_0}\cap \cdots \cap U_{i_p} $
and $j: U_{i_0,\ldots,i_p}\to X$ is the natural embedding. Any $\mathcal H^p$ is a sheaf over $X$ in a  natural way. Consider the following complex 
$$
0\to \mathcal G  \xrightarrow{i} \mathcal H^0  \xrightarrow{d^0} \mathcal H^1\xrightarrow{d^1} \cdots,
$$
where $d^i$ is induced by the \u{C}ech coboundary operator. We have 
$$
H^q(X,\mathcal H^p)\simeq \prod_{i_0,\ldots,i_p} H^q(U_{i_0,\ldots,i_p},\mathcal G).
$$
Therefore, $H^1(X,\mathcal H^0)=\lbrace0 \rbrace$. Let us construct an isomorphism $\tau: H^1(\mathcal U, \mathcal G) \to H^1(X,\mathcal G)$.
Consider the following exact sequence of sheaves over $X$
$$
0\to \mathcal G  \xrightarrow{i}\mathcal H^0  \xrightarrow{d^0} \mathcal Ker d^1\to 0.
$$
and the corresponding long exact sequence
$$
\to H^0(X,\mathcal H^0)\xrightarrow{d^0_*} H^0(X,\mathcal Ker d^1)\xrightarrow{\delta} H^1(X,\mathcal G)\to H^1(X,\mathcal H^0)=\lbrace0 \rbrace.
$$
It follows that $\delta$ induces an isomorphism 
$$
\bar\delta:  H^1(\mathcal U,\mathcal G)= H^0(X,\mathcal Ker d^1)/\operatorname{Im} d^0_*\to H^1(X,\mathcal G).
$$
The proof is complete.$\square$

\medskip

If we have a map of complex-analytic manifolds $p:M\to N$ and $\mathcal R$ is  a sheaf on $M$, then  we denote by $p_*(\mathcal R)$ the direct image of the sheaf $\mathcal R$. 
We will need the following lemma.

\begin{lemma}\label{lem direct images of A, B,T} Let $\pi: \mathcal M\to \mathcal B$ be a superbundle with compact fiber $\mathcal S$ and $\dim (\mathcal B) = n|m$.  Then the sheaves 
	$$
	(\pi_0)_*( \mathcal A^v_{pq}),\quad  (\pi_0)_*( \mathcal C^v_{pq}), \quad (\pi_0)_*( \widehat{\mathcal T}^v_{pq})
	$$ 
	are locally free sheaves of $\mathcal F_{\mathcal B_0}$-modules of rank 
	\begin{align*}
	N_q\dim H^0(\mathcal S_0,(\mathcal A_{\mathcal S})_{p-q}), \quad N_q\dim H^0(\mathcal S_0,(\mathcal C_{\mathcal S})_{p-q}),\quad N_q \dim \mathfrak v(\mathcal S)_{p-q},
	\end{align*}
	respectively. Here $N_q:= \dim \bigwedge^q (m)$, where $\bigwedge (m)$ is the Grassmann algebra with $m$ generators.

\end{lemma}

\noindent {\it Proof.} Assume that the bundle $ \widehat{\mathcal M}$ is trivial over $\mathcal U\subset  \widehat{\mathcal B}$. Let us prove for example that $(\pi_0)_*( \widehat{\mathcal T}^v_{pq})$ is a locally free sheaf of $\mathcal F_{\mathcal B}$-modules. Let $V\subset \mathcal U_0$ be open. By Theorem \ref{theor Kunneth Formula} and Proposition \ref{prop_geometr A,B,T} we have
\begin{align*}
H^0(V, (\pi_0)_*( \widehat{\mathcal T}^v_{pq}))= H^0(\pi_0^{-1}(V), \widehat{\mathcal T}^v_{pq}) \simeq H^0(\pi_0^{-1}(V), \pi_{\mathcal U}^*(\widehat{\mathcal O}_{\mathcal B})_q\otimes \pi_{\mathcal S}^*(\widehat{\mathcal T}_{p-q}))\simeq \\
H^0(V, (\widehat{\mathcal O}_{\mathcal B})_q)\otimes H^0(\mathcal S_0, (\widehat{\mathcal T}_{\mathcal S})_{p-q}).
\end{align*}
We see that $(\pi_0)_*( \widehat{\mathcal T}^v_{pq})|_{\mathcal U_0}$ is a free sheaf of $\mathcal F_{\mathcal B}$-modules. 
The cases $(\pi_0)_*( \mathcal A^v_{pq})$ and $(\pi_0)_*( \mathcal C^v_{pq})$ are similar.$\square$

\section{1st cohomology group with values in $\tilde{\mathcal T}_{p}$}

In this section we will compute cohomology group with values in the sheaves $\mathcal A_{p}$ for all $p$ and $ \mathcal C_{p}$ for $p\ne 2$. More precisely,
we will show that 
$$
H^1(\mathcal M_0, \mathcal A_{p})= \lbrace0\rbrace, \quad p\geq -1, \quad \text{and} \quad H^1(\mathcal M_0, \mathcal C_{p})=\lbrace0\rbrace, \quad p\ne 2.
$$ 
Using these results and Sequence (\ref{exact sequens for T_p})
we will obtain that $H^1(\mathcal M_0, \tilde{\mathcal T}_{p})= \lbrace0\rbrace$, where $p\ne 2$. Further, we will compute
$H^1(\mathcal M_0, \widehat{\mathcal T}^v_{2q})$, $H^1(\mathcal M_0, \widehat{\mathcal T}^h_{2q})$ and $H^1(\mathcal M_0, \widehat{\mathcal T}_{2q})$ and at the end we will show that $H^1(\mathcal M_0, \tilde{\mathcal T}_{2})= \mathbb C$.

\subsection{1st cohomology group with values in $\mathcal A^v_{pq}$,   $\mathcal A^h_{pq}$ and $\mathcal A_{p}$}

As above consider the following Cartan subalgebra in $\mathfrak{gl}_{m|n}(\mathbb C)_{\bar 0}\simeq \mathfrak{gl}_m(\mathbb C)\oplus \mathfrak{gl}_m(\mathbb C)$
$$
\mathfrak{h} := \{\operatorname{diag}(\mu_1,\dots,\mu_{m})\}\oplus
\{\operatorname{diag}(\lambda_1,\dots,\lambda_{n})\}.
$$ 
Recall that the reductive part $R$ of $P$ has the following form see (\ref{eq reductive part R=})
$$
R\simeq \bigoplus_{i=1}^{r+1}\mathfrak{gl}_{k_{i-1}- k_i}(\mathbb{C}) \oplus\bigoplus_{i=1}^{r+1} \mathfrak{gl}_{l_{i-1}- l_i}(\mathbb{C}).
$$ 

We fix the following system of positive roots
$$
\Delta^+=\Delta^+_{1} \cup \Delta^+_{2},
$$
where
$$
\Delta^+_{1}=\{\mu_i-\mu_j, \,\,i<j\} \quad\text{and} \quad \Delta^+_{2}=\{\lambda_p-\lambda_q, \,\,p<q\}.
$$
We denote by $\Phi=\Phi_{1}\cup \Phi_{2}$, where 
$$
\begin{array}{c}
\Phi_{1}= \{\alpha_1,..., \alpha_{m-1}\}, \,\,\,
\alpha_i=\mu_i-\mu_{i+1},\quad \Phi_{2}= \{\beta_1,..., \beta_{n-1}\}, \,\,\,
\beta_j=\lambda_j-\lambda_{j+1},
\end{array}
$$
the system of simple roots. Further, denote by $\mathfrak h^*(\mathbb R)$ the real subspace in $\mathfrak h^*$
spanned by $\mu_j$ and $\lambda_i$. Consider the scalar product $( \,\,,\, )$ in $\mathfrak h^*(\mathbb R)$ such that the vectors  $\mu_j,\lambda_i$ form an orthonormal basis. An element $\delta\in \mathfrak h^*(\mathbb R)$ is called {\it dominant} if $(\delta, \alpha)\ge 0$ for all $\alpha \in \Delta^+$.   Denote by $\zeta$ the half-sum of positive roots. This is 
\begin{equation}\label{eq half-sum of positive roots}
\zeta = \frac{1}{2} \sum_{\alpha\in \Delta^+} \alpha = \frac{1}{2}  \sum_{i=1}^m (m-2i +1) \mu_i + \frac{1}{2}  \sum_{j=1}^n (n-2j +1) \lambda_i. 
\end{equation}
Following \cite{Bott}, we say that $\delta\in \mathfrak h^*(\mathbb R)$ has {\it index $1$} if $(\delta+ \zeta, \alpha) > 0$ for all  $\alpha \in \Delta^+$ except for one root $\beta\in \Delta^+$ for which  $(\delta+\zeta, \beta) < 0$. We call $\delta\in \mathfrak h^*(\mathbb R)$ {\it singular } if $(\delta+ \zeta, \alpha) = 0$ for some  $\alpha \in \Delta^+$.

 Further we will identify $\mathfrak{gl}_{k_{i-1}- k_i}(\mathbb{C})$ or $\mathfrak{gl}_{l_{i-1}- l_i}(\mathbb{C})$ with the corresponding subalgebra in $R$. And we will mark with the superscript $i$ an element in $\{\mu_j\}$ or $\{\lambda_j\}$  if this element is in   
$$
(\mathfrak{h} \cap\mathfrak{gl}_{k_{i-1}- k_i}(\mathbb{C}))^*\quad \text{or} \quad 
(\mathfrak{h} \cap\mathfrak{gl}_{l_{i-1}- l_i}(\mathbb{C}))^*,
$$ 
respectively. 
For example if $\mu_j\in (\mathfrak{h} \cap\mathfrak{gl}_{k_{i-1}- k_i}(\mathbb{C}))^*$ we will write $\mu^i_j$ instead of $\mu_j$.

From Proposition \ref{prop_geometr A,B,T}, Theorem \ref{theor vect fields on flag} and Corollary \ref{cor H^0(A_p)=} we obtain the following result. 

\begin{proposition}
	Assuming (\ref{eq generic type condition for flags}) and (\ref{condition flag type
		1}), we have 
	$$
	(\pi_0)_*( \widehat{\mathcal T}^v_{pq})=0 \quad \text{and} \quad (\pi_0)_*( \mathcal A^v_{pq})=0,
	$$
	for  $p-q\ne -1,0$.
\end{proposition}

\noindent {\it Proof.}  Under assumptions of the proposition by Theorem \ref{theor vect fields on flag}  we have the equality $H^0(\mathcal S_0, (\tilde{\mathcal T}_{\mathcal S})_{p-q})=0$ and by Corollary \ref{cor H^0(A_p)=} we get $H^0(\mathcal S_0, (\mathcal A_{\mathcal S})_{p-q})=0$. The result follows from Proposition \ref{prop_geometr A,B,T}.$\square$

\medskip

As above we denote by $P_{\mathcal B}$ the underlying Lie group of the stabilizer of the origin of $Z_{I_1}$, where $Z_I = (Z_{I_j})$ are as in (\ref{eq I_0 and I_1 of Z_I}). Denote also by $R_{\mathcal B}$ the reductive part of $P_{\mathcal B}$ and by $\phi_{\mathcal B}$ the representation of $P_{\mathcal B}$ corresponding to the homogeneous locally free sheaf $(\mathcal O _{\mathcal B})_1$. In \cite{onigl} it was proved that the representation $\phi_{\mathcal B}$ is completely reducible and  
$$
\phi_{\mathcal B}|_{R_{\mathcal B}}= \rho_1^*\otimes \varsigma_2+ \sigma_1^*\otimes \varrho_2,
$$
 where $\rho_1$ and $\sigma_1$ are standard representations of  $GL_{m-k_1}(\mathbb C)$ and $GL_{n-l_1}(\mathbb C)$, respectively; $\varrho_2$ and $\varsigma_2$ are standard representations of $GL_{k_1}(\mathbb C)$ and $GL_{l_1}(\mathbb C)$, respectively. Note that $\phi_{\mathcal B}$ is equal to $\theta$ from Lemma \ref{lemma_reps of O_11 nd O_10} for $r=1$.

\begin{lemma}\label{lem pi_1(T)} We assume(\ref{eq generic type condition for flags}) and (\ref{condition flag type
		1}). 	The locally free sheaves of $\mathcal F_{\mathcal B}$-modules  
	$$
	(\pi)_*(\widehat{\mathcal T}^v_{-10}),\quad (\pi)_*(\widehat{\mathcal T}^v_{00}),\quad  (\pi)_*(\mathcal A^v_{-10}), \quad  (\pi)_*(\mathcal A^v_{00})
	$$ 
 on $\mathcal B_0$ are homogeneous.  The  corresponding representations of  $P_{\mathcal B}$ are completely reducible and they respectively are:
\begin{align*}
\begin{split}
&\chi_{\widehat{\mathcal T}^v_{-10}}|_{R_{\mathcal B}} =\varrho_2^*\otimes\varsigma_2+  \varrho_2\otimes\varsigma_2^*,\quad \chi_{\widehat{\mathcal T}^v_{00}}|_{R_{\mathcal B}} =
Ad_{\varrho_2}+ Ad_{\varsigma_2}+1+1,\\ 
&\chi_{{\mathcal A}^v_{-10}}|_{R_{\mathcal B}} =
\varrho_2^*\otimes\varsigma_2+  \varrho_2\otimes\varsigma_2^*,\quad \chi_{{\mathcal A}^v_{00}}|_{R_{\mathcal B}} =
1+1,
\end{split}
\end{align*}
where $1$ is the trivial one dimensional representation, $Ad_{\varrho_2}$ and $Ad_{\varsigma_2}$ are adjoint representations of $GL_{k_1}(\mathbb C)$ and $GL_{l_1}(\mathbb C)$, respectively.

\end{lemma}

\medskip

\noindent{\it Proof.} The proposition follows from Proposition \ref{prop_geometr A,B,T}, Theorem \ref{theor vect fields on flag}, the isomorphism $\widehat{\mathcal T}^v_{-10} \simeq  \mathcal A^v_{-10}$ and Corollary \ref{cor H^0(A_p)=}. Further a direct computation shows that the nilpotent part of $P_{\mathcal B}$ acts trivially. Hence we obtain completely reducibility $\Box$ 
\medskip

\begin{corollary} 
	The homogeneous locally free $\mathcal F_{\mathcal B}$-sheaves $(\pi)_*(\widehat{\mathcal T}^v_{pq})$ and $(\pi)_*(\mathcal A^v_{pq})$ correspond to the following representations of $P_{\mathcal B}$
	$$
	\bigwedge^{q}\phi_{\mathcal B}\otimes \chi_{\widehat{\mathcal T}^v_{p-q,0}}, \quad \bigwedge^{q}\phi_{\mathcal B}\otimes \chi_{{\mathcal A}^v_{p-q,0}},
	$$ 
	respectively.
\end{corollary}

\medskip

\noindent{\it Proof.} The result is a consequence of the following observation
$$
\widehat{\mathcal T}^v_{pq} \simeq \bigwedge^{q}  \widehat\pi^* ((\widehat{\mathcal O}_{\mathcal B})_1) \otimes \widehat{\mathcal T}^v_{p-q,0}, \quad {\mathcal A}^v_{pq} \simeq \bigwedge^{q}  \widehat\pi^* ((\widehat{\mathcal O}_{\mathcal B})_1) \otimes {\mathcal A}^v_{p-q,0}. \Box
$$

In \cite{V_funk} the following Lemma was proved. 

\begin{lemma}\label{lem H^0(flag)=0}
Assuming (\ref{eq generic type condition for flags}), we have 	\begin{align*}
	H^0(\mathcal M_0, \tilde{\mathcal O}_p)=
	\left\lbrace
	\begin{array}{c}
	\lbrace0\rbrace, \,\,\, p\ne0,\\
	\mathbb C, \,\,\, p=0.
	\end{array}\right.
	\end{align*}
\end{lemma}

Now we need the following theorem.

\begin{theorem}\label{theor H^1(O)=0}
	We assume(\ref{eq generic type condition for flags}) and (\ref{condition flag type
		1}). Then 
$$
H^1(\mathcal M_0,\mathcal O)=\lbrace0\rbrace \quad \text{and} \quad
	H^1(\mathcal M_0,\tilde{\mathcal O})=\lbrace0\rbrace.
$$
\end{theorem}

\noindent {\it Proof.} We use induction on $r$. For the case  $\mathcal M = \mathbf{Gr}_{m\vert n,k\vert l}$ by \cite{onigl} we have $\tilde{\mathcal O}\simeq \bigwedge\mathcal E_{\phi_{\mathcal B}}$, where $\mathcal E_{\phi_{\mathcal B}}$ is the homogeneous locally free sheaf corresponding to the following representation of $P_{\mathcal B} $
$$
\phi_{\mathcal B}|_{R_{\mathcal B}}= \rho_1^*\otimes \varsigma_2+ \sigma_1^*\otimes \varrho_2.
$$ 
Here $\rho_1$, $\varrho_2$, $\sigma_1$, $\varsigma_2$ are standard representations of  $GL_{m-k}(\mathbb C)$, $GL_{k}(\mathbb C)$, $GL_{n-l}(\mathbb C)$, $GL_{l}(\mathbb C)$, respectively. Let us show that $H^1(\mathcal M_0,\bigwedge^q\mathcal E_{\phi_{\mathcal B}})=\lbrace0\rbrace$ for $q\geq 0$.

Any weight of the representation  $\bigwedge^q \phi_{\mathcal B}$ has the form $\gamma= \gamma_0+\gamma_1$, where 
$$
\gamma_0= -\mu_{i_1}^1-\cdots -\mu_{i_a}^1+\mu_{j_1}^2+\cdots +\mu_{j_b}^2, \quad \gamma_1= -\lambda_{i_1}^1-\cdots -\lambda_{i_b}^1+\lambda_{j_1}^2+\cdots +\lambda_{j_a}^2,\quad a+b= q.
$$
Assume that $\gamma_0$ is dominant. Then under assumption $0<k<m$, we get that $a=b=0$. In this case $\gamma_1$ cannot have index $1$. Similarly we get that if $\gamma_1$ is dominant $\gamma_0$ cannot have index $1$. Summing up under our assumptions the weight $\gamma$ cannot have index $1$ and therefore  by Borel-Weil-Bott theorem we have
$$
H^1(\mathcal M_0,\tilde{\mathcal O}_q) =
H^1(\mathcal M_0,\bigwedge^q\mathcal E_{\phi_{\mathcal B}})=\lbrace0\rbrace, \quad q\geq 0.
$$

Further we have the following exact sequence 
$$
0\to \mathcal O_{(p+1)}\to \mathcal O_{(p)} \to \tilde{\mathcal O}_{p}\to 0.
$$
For enough big $p$ we have the natural isomorphism $\tilde{\mathcal O}_p\simeq \mathcal O_{(p)}$ and therefore for enough big $p$ we have $ H^1(\mathcal M_0,\mathcal O_{(p)} ) \simeq H^1(\mathcal M_0,\tilde{\mathcal O}_{p} )=\lbrace0\rbrace$. By induction we obtain $H^1(\mathcal M_0,\mathcal O_{(p)})=\lbrace0\rbrace$  for any $p$.

Now assume that $\mathcal M$ is a flag supermanifold of length  $r>1$. Then $\mathcal M$ is a super\-bundle with base space  $\mathcal B= \mathbf{Gr}_{m\vert n,k_1\vert l_1}$ and with fiber $\mathcal S$, where $\mathcal S$ is  a flag supermanifold  of length $r-1$. Denote by $\tilde{\pi}: \tilde{\mathcal M} \to \tilde{\mathcal B}$ the  projection of the split superbundle. 
Let $\{\mathcal U \}$ be a covering of $\tilde{\mathcal B}$ and assume that the bundle $\tilde{\mathcal M}$ is trivial over any $\mathcal U$. Assume by induction that $H^1(\mathcal S_0, \tilde{\mathcal O}_{\mathcal S})= \lbrace0\rbrace$.  Then by Theorem \ref{theor Kunneth Formula} and by induction we have
$$
H^1(\pi^{-1}(\mathcal U),\tilde{\mathcal O}\vert_{\pi^{-1}(\mathcal U)}) = H^0(\mathcal U_0, \tilde{\mathcal O}_{\mathcal B} \vert_{\mathcal U_0})\otimes H^1(\mathcal S_0, \tilde{\mathcal O}_{\mathcal S})= \lbrace0\rbrace.
$$ 
Hence by Lemma \ref{lem H^1(covering) = H^1} 
$$
H^1(\mathcal M_0,\tilde{\mathcal O})\simeq H^1(\lbrace \pi^{-1}( \mathcal U) \rbrace,\tilde{\mathcal O})\simeq H^1(\lbrace U\rbrace,\pi_*(\tilde{\mathcal O})).
$$ 

By Lemma \ref{lem H^0(flag)=0} we have $H^0(\mathcal S_0, \tilde{\mathcal O}_{\mathcal S})\simeq \mathbb C$. Therefore by Theorem \ref{theor Kunneth Formula}, $\tilde{\pi}_*(\tilde{\mathcal O})= \tilde{\mathcal O}_{\mathcal B}$. Hence 
$$
H^1(\mathcal M_0,\tilde{\mathcal O})= H^1(\mathcal B_0,\tilde{\mathcal O}_{\mathcal B})=\lbrace0\rbrace.
$$ 
The argument in the  case of $\mathcal O$ is similar to the case $r=1$.$\square$

\medskip

In \cite[Propositions 1,2,3]{onirig} the following lemma was proved.

\begin{lemma}\label{lem Oni H^1(A_p)=0}
	Assume that $r=1$, and (\ref{eq generic type condition for flags}) and (\ref{condition flag type
		1}) hold true. Then for $\mathcal M=\mathbf{Gr}_{m\vert n, k\vert l}$ we have
	$$
	H^1(\mathcal M_0, \mathcal A_{p})=\lbrace0\rbrace, \quad p\geq -1.
	$$
\end{lemma}

Now we can prove the following theorem.

\begin{theorem}\label{theor A_p= and A_pq =} 
	Assume that $r\geq 2$ and that (\ref{eq generic type condition for flags}) and (\ref{condition flag type 1}) hold true.
	Then
	\begin{align*}
	H^1(\mathcal M_0,\mathcal A_{pq}^v)&=
	\left \lbrace
	\begin{array}{l}
	\mathbb C^2, \,\,\, (p,q)=(0,1),  \\
	\lbrace 0\rbrace , \,\,\, (p,q)\ne (0,1),
	\end{array}
	\right.\\
	H^1(\mathcal M_0,\mathcal A_{pq}^h)&= \lbrace 0\rbrace \,\,\,\text{for all }\,\,\, p,q.\\
		H^1(\mathcal M_0, \mathcal A_{p})&=\lbrace0\rbrace\,\,\,\text{for all }\,\,\, p.
	\end{align*}
\end{theorem}

\noindent {\it Proof.} We use induction on $r$. For $r=1$ the sheaves $\mathcal A_{pq}^v$ and $\mathcal A_{pq}^h$ are no defined and $H^1(\mathcal M_0, \mathcal A_{p})=\lbrace0\rbrace$ by Lemma \ref{lem Oni H^1(A_p)=0}. Assume by induction that for a flag supermanifold of the length $r-1$ we have $H^1(\mathcal M_0, \mathcal A_{p})=\lbrace0\rbrace\,\,\,\text{for all }\,\,\, p$. Let us prove the statement for the length $r$. We have the following exact sequence
\begin{equation}\label{eq exact sequence A_{p(q+1)} to}
0\to \mathcal A_{p(q+1)}\to \mathcal A_{p(q)}\to \mathcal A_{pq}\to 0,
\end{equation}
where $\mathcal A_{pq}= \mathcal A_{pq}^v+ \mathcal A_{pq}^h$. Let us compute $1$-cohomology group for the sheaves $\mathcal A_{pq}^v$ and $\mathcal A_{pq}^h$.

 {\bf Step 1: $1$st cohomology group with values in  $\mathcal A_{pq}^h$.} Let $\{\mathcal U \}$ be a covering of $\mathcal B$ and assume that the bundle $\mathcal M$ is trivial over any $\mathcal U$. By Theorem \ref{theor H^1(O)=0} we have
$H^1(\mathcal S_0, (\tilde{\mathcal O}_{\mathcal S})_{p-q})=\lbrace0\rbrace$ for $p-q\geq 0$. Therefore, by Theorem \ref{theor Kunneth Formula} and by Proposition \ref{prop_geometr A,B,T} we have $H^1(\pi^{-1}(\mathcal U_0), \mathcal A_{pq}^h)=0$ and hence by Lemma \ref{lem H^1(covering) = H^1} we get 
$$
H^1(\mathcal M_0,\mathcal A_{pq}^h)=H^1(\{\pi^{-1}(\mathcal U_0 )\}, \mathcal A_{pq}^h).
$$ 
Further by Theorem \ref{theor Kunneth Formula}, by Proposition \ref{prop_geometr A,B,T} and by Lemma \ref{lem H^0(flag)=0} we get
\begin{align*}
H^0 (\pi^{-1}(\mathcal U_0), \mathcal A_{pq}^h) =  H^0(\mathcal S_0, \bigwedge^{p-q}(\tilde{\mathcal{O}}_{\mathcal B})_{1})\otimes H^0(\mathcal U_0, (\mathcal A_{\mathcal{B}})_{q} )  = \{0\} \quad \text{for} \,\,  p\ne q,
\end{align*}
and 
\begin{align*}
H^0 (\pi^{-1}(\mathcal U_0), \mathcal A_{pp}^h) = H^0(\mathcal U_0, (\mathcal A_{\mathcal{B}})_{p} ) \quad \text{for} \,\,  p=q.
\end{align*}
Therefore, $H^1 (\{\pi^{-1}{(\mathcal U_0)} \}, \mathcal A_{pp}^h) = H^1(\{\mathcal U\}, (\mathcal A_{\mathcal{B}})_{p})$ and hence $H^1 (\mathcal M_0, \mathcal A_{pp}^h) =  \{ 0\}  $ by Lemma \ref{lem Oni H^1(A_p)=0}. Summing up we proved that 
$$
H^1 (\mathcal M_0, \mathcal A_{pq}^h) =  \{ 0\}, \quad p,q \geq -1. 
$$

{\bf Step 2: $1$st cohomology group with values in  $\mathcal A_{pq}^v$.}
Further by the induction assumption, $H^1(\mathcal S_0, (\mathcal A_{\mathcal S})_{p-q})=\lbrace0\rbrace$. Hence  by Proposition \ref{prop_geometr A,B,T} and by Theorem \ref{theor Kunneth Formula} we have $H^1(\pi^{-1}(\mathcal U_0), \mathcal A_{pq}^v)=\lbrace0\rbrace$. Therefore by Lemma \ref{lem H^1(covering) = H^1} to compute $1$st cohomology we can use the covering $\{ \pi^{-1}(\mathcal U_0)\}$. Now consider the homogeneous locally free sheaf $\pi_*(\mathcal A_{pq}^v)$ of $\mathcal F_{\mathcal B}$-modules. The corresponding representation of $P_{\mathcal B}$ is $\bigwedge^{q}\phi_{\mathcal B}\otimes \chi_{{\mathcal A}^v_{p-q,0}}$, see Lemma \ref{lem pi_1(T)} and its Corollary. Let us apply the Borel-Weil-Bott Theorem.

Consider first of all the case $p-q=-1$. Any weight of $\bigwedge^{q}\phi_{\mathcal B}\otimes \chi_{{\mathcal A}^v_{-1,0}}$  has the form $\Lambda= \Lambda_0+\Lambda_1$, where
$$
\begin{array}{l}
\Lambda_0= -\mu_{i_1}^1-\cdots-\mu_{i_a}^1+\mu_{j_1}^{s}+\cdots+\mu_{j_b}^{s}+
\left \lbrace
\begin{array}{l}
(1.)\,\,\, -\mu_{i}^{s},\\
(2.)\,\,\, \mu_{j}^{s},\\
\end{array}
\right. \\
\Lambda_1= -\lambda_{j_1}^1-\cdots-
\lambda_{i_b}^1+\lambda_{i_1}^{s}+
\cdots+\lambda_{i_a}^{s}+
\left \lbrace
\begin{array}{l}
(1.)\,\,\, \lambda_{j}^{s},\\
(2.)\,\,\, -\lambda_{i}^{s},\\
\end{array}
\right.
\end{array} \,\,a+b=q,\,\, s>1.
$$
Assume that $\Lambda$ has index $1$. Consider the case when $\Lambda_0$ is dominant and   $\Lambda_1$ has index $1$. In case  $(1)$ it follows that either $(1.A)$ $a=0$, $b=1$, $\Lambda_0=0$; or $(1.B)$ $a=b=0$. In case $(1.B)$ the weight $\Lambda_1$ is singular. In case $(1.A)$ the weight $\Lambda_1$ is singular except of  $\Lambda_1=  -\lambda_{n-l_1}^1+\lambda_{n-l_1+1}^2$, which has index $1$. In case $(2)$,  $\Lambda_0$ is not dominant. Summing up,  $\Lambda= -\lambda_{n-l_1}^1+\lambda_{n-l_1+1}^2$ is the unique weight of index $1$. We see that $\Lambda$ is a highest weight of the representation  $\bigwedge^{q}\phi_{\mathcal B}\otimes \chi_{{\mathcal A}^v_{-1,0}}$ for $q=1$. In case if $\Lambda_1$ is dominant and  $\Lambda_0$ has index $1$ similarly as above we get another highest weight  $\Lambda= -\mu_{m-k_1}^1+\mu_{m-k_1+1}^2$ of $\bigwedge^{q}\phi_{\mathcal B}\otimes \chi_{{\mathcal A}^v_{-1,0}}$ of index $1$ for $q=1$. By the Borel-Weil-Bott Theorem we get
\begin{align*}
H^1(\mathcal M_0,\mathcal A_{q-1,q}^v) =  H^1(\mathcal B_0,\pi_*(\mathcal A_{q-1,q}^v))=  \left \lbrace
\begin{array}{l}
 \mathbb C^2,\,\, q = 1;\\
 \lbrace0\rbrace,\,\, q\ne 1.\\
\end{array}
\right.
\end{align*}

Consider now the case $p-q=0$. Any weight of the representation $\bigwedge^{q}\phi_{\mathcal B}\otimes \chi_{{\mathcal A}^v_{0,0}}$ corresponding to the locally free sheaf $\pi_*(\mathcal A_{qq}^v)$  has the form $\Lambda= \Lambda_0+\Lambda_1$, where
$$
\begin{array}{l}
\Lambda_0= -\mu_{i_1}^1-\cdots-\mu_{i_a}^1+\mu_{j_1}^{s}+\cdots+\mu_{j_b}^{s},
\\
\Lambda_1= -\lambda_{j_1}^1-\cdots-
\lambda_{j_b}^1+\lambda_{i_1}^{s}+
\cdots+\lambda_{i_a}^{s},
\end{array} \,\,a+b=q, \,\, s>1.
$$
Assume that $\Lambda$ has index $1$. Consider the case when  $\Lambda_0$ is dominant and  $\Lambda_1$ has index $1$. Hence $a=b=0$ and  $\Lambda_1=0$. We get a contradiction with the assumption that $\Lambda$ has index $1$. The case when  $\Lambda_1$ is dominant and  $\Lambda_0$ has index $1$ is similar.  Therefore the representation $\bigwedge^{q}\phi_{\mathcal B}\otimes \chi_{{\mathcal A}^v_{0,0}}$ does not have weights of index $1$. We get that  
$$
H^1(\mathcal M_0,\mathcal A_{qq}^v) = H^1(\mathcal B_0, \pi_*(\mathcal A_{qq}^v))= \lbrace 0\rbrace.
$$
For  $(p,q)$ such that $p-q\ne -1,0$ we get using Corollary \ref{cor H^0(A_p)=} 
$$
H^0(\pi^{-1}(\mathcal U), \mathcal A_{pq}^v)=\lbrace0\rbrace.
$$
Hence,
$$
H^1(\mathcal M_0,\mathcal A_{pq}^v) =  \lbrace 0\rbrace. 
$$

{\bf Step 3: $1$st cohomology group with values in  $\mathcal A_{p}$.}
Now from the exact sequence (\ref{eq exact sequence A_{p(q+1)} to})
we get that  $H^1(\mathcal M_0, \mathcal A_{p})= \lbrace0\rbrace$, where $p\ne 0$. Further, 
$$
H^1(\mathcal M_0,\mathcal A_{0(1)})= H^1(\mathcal M_0,\mathcal A_{01})=  H^1(\mathcal M_0,\mathcal A_{01}^v \oplus \mathcal A_{01}^h)= \mathbb C^2.
$$ 
We have the following long exact sequence
$$
\begin{array}{ll}
0\to & H^0(\mathcal M_0,\mathcal A_{0(1)}) \to H^0(\mathcal M_0,\mathcal A_{0(0)})\to H^0(\mathcal M_0,\mathcal A_{00})\to\\
&H^1(\mathcal M_0,\mathcal A_{0(1)}) \to H^1(\mathcal M_0,\mathcal A_{0(0)})\to H^1(\mathcal M_0,\mathcal A_{00})=0.
\end{array}
$$

From Corollary of Theorem \ref{theor vect fields on flag} we know that $H^1(\mathcal M_0,\mathcal A_{0})= \mathbb C^2$. A direct calculation shows that  
$$
H^0(\mathcal M_0,\mathcal A_{0(1)})= \lbrace 0\rbrace,\quad H^0(\mathcal M_0,\mathcal A_{0(0)})= \mathbb C^2.
$$ 
Let us compute $H^0(\mathcal M_0,\mathcal A_{00})$. Clearly, $\mathcal A_{00}= \widehat{\mathcal O}_{11}^*\otimes \widehat{\mathcal O}_{11}+ \widehat{\mathcal O}_{10}^*\otimes \widehat{\mathcal O}_{10}$. As above using the Borel-Weil-Bott Theorem, we obtain  
$$
H^0(\mathcal M_0,\widehat{\mathcal O}_{11}^*\otimes \widehat{\mathcal O}_{11})= H^0(\mathcal B_0,\tilde{\mathcal O}_{\mathcal B}^*\otimes \tilde{\mathcal O}_{\mathcal B})= \mathbb C^2, \quad  H^0(\mathcal M_0,\widehat{\mathcal O}_{10}^*\otimes \widehat{\mathcal O}_{10}) = H^0(\mathcal B_0, \widetilde{\mathcal W}_0)=\mathbb C^2.
$$
%(Compare with the proof of Theorem \ref{theor vect fields on flag}.)
Hence, $H^1(\mathcal M_0,\mathcal A_{0(0)})= \lbrace 0\rbrace$. From the exact sequence
\begin{align*}
\lbrace0 \rbrace= H^1(\mathcal M_0,\mathcal A_{0(0)}) \to H^1(\mathcal M_0,\mathcal A_{0(-1)})\to H^1(\mathcal M_0,\mathcal A_{0-1})=\lbrace0 \rbrace
\end{align*}
it follows that $H^1(\mathcal M_0,\mathcal A_{0})=H^1(\mathcal M_0,\mathcal A_{0(-1)})=\lbrace0 \rbrace$. The proof is complete.$\square$

\medskip

\subsection{1st cohomology group with values in $\mathcal C^v_{pq}$,   $\mathcal C^h_{pq}$ and $\mathcal C_{p}$}

By Theorem \ref{theor vect fields on flag}, Corollary \ref{cor H^0(A_p)=} and   Theorem \ref{theor A_p= and A_pq =} we get.

\begin{lemma}\label{lem H^0(B_p)=}
Assume that $r\geq 1$ and that (\ref{eq generic type condition for flags}) and (\ref{condition flag type 1}) hold true. Then
	$$
	H^0(\mathcal M_0, \mathcal C_{p})=
	\left \lbrace
	\begin{array}{l}
	 \mathfrak {sl}_m(\mathbb C)\oplus\mathfrak {sl}_n(\mathbb C),\,\,\, \,\,\, p=0,  \\
	\lbrace 0\rbrace , \,\,\, p\ne 0.
	\end{array}
	\right.
	$$
\end{lemma}
In case $r=1$ this result was obtained in \cite{onigl}. 

\medskip

\noindent{\it Proof.} The result is a consequence of Theorem \ref{theor vect fields on flag}, Corollary \ref{cor H^0(A_p)=}, Theorem \ref{theor A_p= and A_pq =} and the following exact sequence
$$
0\to H^0(\mathcal M_0, \mathcal A_p) \to H^0(\mathcal M_0, \tilde{\mathcal T}_p) \to H^0(\mathcal M_0, \mathcal C_p) \to H^1(\mathcal M_0, \mathcal A_p) =0.\Box
$$

\medskip

By Proposition \ref{prop_geometr A,B,T}, Theorem \ref{theor Kunneth Formula} and by Lemma \ref{lem H^0(B_p)=}, we get. 

\begin{lemma}\label{lem pi_*(C_pq)=0, p ne q}
	We have $\pi_*(\mathcal C_{pq}^v)=\lbrace0\rbrace$ for $p\ne q$.
\end{lemma}

 Further we need the following lemma, which is a consequence of Lemma \ref{lem H^0(B_p)=} and  Proposition \ref{prop_geometr A,B,T}. 

\begin{lemma}\label{lem representations of B_00^v}
The representations of $P_{\mathcal B}$ corresponding to the sheaves $\pi_*(\mathcal C_{p0}^v)$ and $\pi_*(\mathcal C_{pq}^v)$ have the following form
$$
\chi_{\mathcal C_{p0}^v}\vert_R = \left \lbrace
\begin{array}{l}
 0,\,\,\, \,\,\, p\ne 0,  \\
Ad_{\varrho_2}+ Ad_{\varsigma_2}, \,\,\, p= 0.
\end{array}
\right., 
\quad \chi_{\mathcal C_{pq}^v}\vert_R= \bigwedge^q\phi_{\mathcal B}\otimes \chi_{\mathcal C_{p-q,0}^v},
$$
respectively.
\end{lemma}

In \cite{onirig} for $\mathbf{Gr}_{m\vert n, k\vert l}$ the following lemma was proved. 

\begin{lemma}\label{lem H^1(B_p)=}
	Assume that $r=1$, $0<k<m$, $0<l<n$ and $(k,l)\ne (1,n-1), (m-1,1)$. Then 
	$$
	H^1(\mathcal M_0, \mathcal C_{p})=
	\left \lbrace
	\begin{array}{l}
	\mathbb C^2,\,\,\, \,\,\, p=2,  \\
	\lbrace 0\rbrace , \,\,\, p\ne 2.
	\end{array}
	\right.
	$$
\end{lemma}

We are ready to prove the following theorem. 

\begin{theorem}\label{theor H^1(B^h_pq)=}
	Assume that $r\geq 2$ and that (\ref{eq generic type condition for flags}) and (\ref{condition flag type 1}) hold true. Then 
	$$
	H^1(\mathcal M_0, \mathcal C_{p})=\lbrace0\rbrace \quad \text{for} \quad p\ne 2
	$$ 
	 and
	$$
	H^1(\mathcal M_0,\mathcal C_{pq}^h)=
	\left \lbrace
	\begin{array}{l}
	\lbrace0\rbrace\,\,\, (p,q)\ne (2,2),\\
	\mathbb C^2, \,\,\, (p,q)= (2,2),\\
	\end{array}
	\right.\quad H^1(\mathcal M_0, \mathcal C_{pq}^v)=\lbrace0\rbrace,\quad (p,q)\ne (2,0).
	$$
\end{theorem}

\noindent {\it Proof.} Consider the following exact sequence
\begin{equation}\label{eq exact B^v_p to B_p}
0\to \mathcal C_{p}^v\to \mathcal C_{p}\to \mathcal C_{p}^h\to 0.
\end{equation}
Let us compute the 1-st cohomology of the sheaves 
 $\mathcal C_{p}^v$ and $\mathcal C_{p}^h$, where $p\ne 2$.

Recall that by Proposition \ref{prop_geometr A,B,T} we have 
$$
\mathcal C_{pq}^h\simeq  \bigwedge\limits^{p- q}\widehat{\mathcal{O}}_{10}\otimes {\pi}^*((\mathcal C_{\mathcal{B}})_{q}),\quad \mathcal C_{pq}^v\vert_{\widehat\pi^{-1}(\mathcal{U})} \simeq \bigwedge\limits^{q}\widehat{\mathcal O}_{11}\vert_{\mathcal{U}}\otimes \pi_{\mathcal S}^* (\mathcal
C_{\mathcal{S}})_{p-q}.
$$

{\bf Step 1, $1$st cohomology group with values in $\mathcal C_{pq}^h$.} By Theorem \ref{theor H^1(O)=0} and by Lemma \ref{lem H^1(covering) = H^1} to compute $H^1(\mathcal M_0,\mathcal C_{pq}^h )$ we can use the covering  $\lbrace\pi^{-1} (\mathcal U)\rbrace$. Further, by Lemma \ref{lem H^0(flag)=0} and Theorem \ref{theor Kunneth Formula} we have
$H^0(\lbrace \pi^{-1} (\mathcal U)\rbrace,\mathcal C_{pq}^h )= \lbrace0\rbrace$ for $p\ne q$ and $H^0(\lbrace \pi^{-1} (\mathcal U)\rbrace,\mathcal C_{pp}^h )= H^0(\mathcal B_0,(\mathcal C_{\mathcal B})_p )$. Therefore by Lemma \ref{lem H^1(B_p)=} we have
$$
H^1(\mathcal M_0,\mathcal C_{pq}^h)=
\left \lbrace
\begin{array}{l}
\lbrace0\rbrace\,\,\, (p,q)\ne (2,2),\\
\mathbb C^2, \,\,\, (p,q)= (2,2).\\
\end{array}
\right.
$$

{\bf Step 2, $1$st cohomology group with values in $\mathcal C_{pq}^v$ for $p-q\ne 2$.}
For $p-q\ne 2$ by induction we assume that $H^1(\mathcal S_0,(\mathcal C_{\mathcal S})_{p-q})= \lbrace0\rbrace$. Hence in this case we can use Theorem \ref{theor Kunneth Formula} and Lemma \ref{lem H^1(covering) = H^1}. Consider the locally free sheaf $\pi_*(\mathcal C_{pq}^v)$ of $\mathcal F_{\mathcal B}$-modules.  Note that by Lemma \ref{lem H^0(B_p)=} we have $H^0(\mathcal U_0,\pi_*(\mathcal C_{pq}^v))= \lbrace0\rbrace$ for $p-q\ne 0$. Hence, 
$$
H^1(\mathcal M_0,\mathcal C_{pq}^v)= \lbrace0\rbrace \quad \text{for}\quad p-q\ne 2, 0.
$$ 
Further the representation $\chi_{\mathcal C_{pq}^v}$ of $P_{\mathcal B}$ in a fiber of $\pi_*(\mathcal C_{pq}^v)$ is completely reducible, since the nilradical of  $P_{\mathcal B}$ acts trivially. We use the Borel--Weil--Bott Theorem. By Lemma \ref{lem representations of B_00^v} any weight of $\chi_{\mathcal C_{pq}^v}$ for $p-q=0$ has the form $\Lambda= \Lambda_0+\Lambda_1$, where
$$
\begin{array}{l}
\Lambda_0= -\mu_{i_1}^1-\cdots-\mu_{i_a}^1+\mu_{j_1}^2+\cdots+\mu_{j_b}^2+
\left [
\begin{array}{l}
1.\,\,\, \mu_{j}^2-\mu_{i}^2,\\
2.\,\,\, 0,\\
\end{array}
\right. \\
\Lambda_1= -\lambda_{j_1}^1-\cdots-
\lambda_{i_b}^1+\lambda_{i_1}^2+
\cdots+\lambda_{i_a}^2+
\left [
\begin{array}{l}
1.\,\,\, 0,\\
2.\,\,\, \lambda_{j}^2-\lambda_{i}^2,\\
\end{array}
\right.
\end{array} \,\,a+b=q.
$$
The weights $\Lambda_0$ and $\Lambda_1$ are not dominant for any $a,b$. Therefore $\Lambda$ cannot have index $1$ and
$$
H^1(\mathcal M,\mathcal C_{pq}^v)= \lbrace0\rbrace, \,\,\, p-q\ne 2.
$$

{\bf Step 3, $1$st cohomology group with values in $\mathcal C_{pq}^v$ for $p-q= 2$ and with values in $\mathcal C_{p}$.} Note that in this case $q \geq 0$. Assume that $q>0$. Then by Lemma \ref{lem q>0} below the representation of $P$ corresponding to $\mathcal C_{pq}^v$ does not have weights of index $1$. Therefore we have only one possibility $(p,q)=(2,0)$. Summing up, we get the following result 
$$
H^1(\mathcal M_0, \mathcal C_{pq}^v)=\lbrace0\rbrace\quad \text{for} \quad 	 (p,q)\ne (2,0).
$$ 
This implies that $H^1(\mathcal M_0, \mathcal C_{p}^v)=\lbrace0\rbrace$ for $p\ne 2$. Hence we get 
$$
H^1(\mathcal M_0, \mathcal C_{p})=\lbrace0\rbrace\quad \text{for} \quad p\ne 2.\square
$$

\medskip

\begin{lemma}\label{lem q>0} Assume that $r\geq 2$, that (\ref{eq generic type condition for flags}) and (\ref{condition flag type 1}) hold true, $q>0$ and $p-q=2$. Then the representation 
	$\tau^v\otimes \bigwedge^{q}\varphi\otimes \bigwedge^{p-q}\phi$ corresponding to the sheaf $\mathcal C_{pq}^v$ does not have weights of index $1$.	
\end{lemma}

\noindent{\it Proof.} We use notations of Lemma \ref{lemma_reps of O_11 nd O_10} and  Lemma \ref{lemma_reps of T,T^v and T/T^v}. 
Recall that we denoted by $\varphi$ the representation of $P$ corresponding to the homogeneous locally free sheaf $\widehat{\mathcal O}_{11}$ and by $\psi$ the representation of $P$ corresponding to the homogeneous locally free sheaf $\widehat{\mathcal O}_{10}$. A weight $\Lambda$ of the representation   
$$
\tau^v\otimes \bigwedge^{q}\varphi\otimes
\bigwedge^{2}\psi
$$
of $P$ has the form $\Lambda= \Lambda_0 + \Lambda_1$, where
$$
\Lambda_0= -\mu_{i_1}^1-\cdots-\mu_{i_a}^1+\mu_{j_1}^{\geq 2}+\cdots+\mu_{j_b}^{\geq 2}+ 
\left\{ 
\begin{array}{ll}
-\mu_{i}^{\nu_1}-\mu_{j}^{\nu_2}\,& (A),\\
-\mu_{i}^{\nu_1} + \mu_{j}^{\geq 3}\,& (B)\\
+\mu_{i}^{\geq 3} + \mu_{j}^{\geq 3}\,& (C)\\
\end{array}\right.
+
\left\{ 
\begin{array}{ll}
\mu_{s}^{\nu_3}-\mu_{t}^{\geq 3}\,& (1),\\
0\,& (2)
\end{array}\right.,
$$
and 
$$
\Lambda_1= -\lambda_{j_1}^1-\cdots-
\lambda_{i_b}^1+\lambda_{i_1}^{\geq 2}+
\cdots+\lambda_{i_a}^{\geq 2}+  
\left\{ 
\begin{array}{ll}
\lambda_{i}^{\geq 3}+\lambda_{j}^{\geq 3}\,& (A),\\
-\lambda_{i}^{\kappa_1} + \lambda_{j}^{\geq 3}\,& (B)\\
-\lambda_{i}^{\kappa_1} - \lambda_{j}^{\kappa_2}\,& (C)\\
\end{array}\right.
+
\left\{ 
\begin{array}{ll}
0\,& (1),\\
\lambda_{i}^{\kappa_3}-\lambda_{j}^{\geq 3}\,& (2)
\end{array}\right..
$$
Here $a+b=q$ and $1< \nu_i, \kappa_i< r+1$. 

Assume that $\Lambda_0$ is dominant and $\Lambda_1$ has index $1$. Since $q>0$ we have $a>0$ or $b>0$. Let us write $\Lambda_0$ in the following form
$$
\Lambda_0 = \sum \alpha_i \mu_i.
$$
Then $\Lambda_0$ is dominant if and only if $\alpha_i \geq \alpha_{i+1}$ for any $i$. If $a>0$ and $\Lambda_0$ is dominant we have only one possibility: $(A.1)$ with $a=1$, $b=0$, $r=2$. However in this case $\Lambda_1$ is singular.

Now consider the case $a=0$ and $b>0$. In case $(A.1)$ if $\Lambda_0$ is dominant we have the following possibilities:
\begin{enumerate}
	\item[(I)] $b=1$ and $\Lambda_0 = -\mu^{r+1}_m$; 
	\item[(II)] $b=2$ and $\Lambda_0 = 0$.
\end{enumerate}
In case (I) the weight $\Lambda_1$ is singular. In case (II) the weight  $\Lambda_1$ is singular or has index greater then $1$. 

In case $(A.2)$ if $\Lambda_0$ is dominant we have the following possibility:
\begin{enumerate}
	\item[(I)] $b=2$ and $\Lambda_0 = 0$.
\end{enumerate}
In this case the weight $\Lambda_1$ is singular or has index greater then $1$. 

In other cases, this is the cases $(B.1)$, $(B.2)$, $(C.1)$ and $(C.2)$,  the weight $\Lambda_0$ is not dominant if $b>0$. The proof is complete.$\Box$

From Theorem \ref{theor A_p= and A_pq =} and Theorem \ref{theor H^1(B^h_pq)=} the following result follows.

\begin{theorem}\label{theor H^1(T^h_pq,T^v_pq)}
	Assume that $r\geq 2$ and that (\ref{eq generic type condition for flags}) and (\ref{condition flag type 1}) hold true. Then 
	$$
	H^1(\mathcal M_0, \mathcal T_{p})=\lbrace0\rbrace,\quad \text{ where} \quad p\ne 2,
	$$ 
	and
	\begin{align*}
	&H^1(\mathcal M_0,\widehat{\mathcal T}_{pq}^h)=
	\left \lbrace
	\begin{array}{l}
	\lbrace0\rbrace\,\,\, (p,q)\ne (2,2),\\
	?, \,\,\, (p,q)= (2,2),\\
	\end{array}
	\right.\\
	&H^1(\mathcal M_0, \widehat{\mathcal T}_{pq}^v)=
	\left \lbrace
	\begin{array}{l}
	\lbrace0\rbrace\,\,\, (p,q)\ne (2,0), \,(0,1),\\
	\mathbb C^2, \,\,\, (p,q)= (0,1),\\
	?, \,\,\, (p,q)= (2,0).
	\end{array}
	\right.
	\end{align*}
\end{theorem}

\medskip

\noindent{\it Proof.} The result follows from the exact sequence (\ref{exact sequens for T_p}), Theorem \ref{theor A_p= and A_pq =}, Theorem \ref{theor H^1(B^h_pq)=},  Lemma \ref{lem pi_*(C_pq)=0, p ne q} and Lemma \ref{Lemma_diag for A,B}.

\subsection{ $1$st cohomology group with values in the sheaf  $\widehat{\mathcal T}_{22}$}

Let us prove first the following theorem.

\begin{theorem}\label{theor H^1(T^h_22)=C}
Assuming (\ref{eq generic type condition for flags}) and (\ref{condition flag type 1}), we have	
	$H^1(\mathcal M_0, \widehat{\mathcal T}_{22}^h)=\mathbb C$.
\end{theorem}

\smallskip

\noindent {\it Proof.} By Proposition \ref{prop_geometr A,B,T} we have
$$
\widehat{\mathcal T}^h_{pq} \simeq \bigwedge\limits^{p- q}\widehat{\mathcal{O}}_{10}\otimes {\pi}^*((\widehat{\mathcal T}_{\mathcal{B}})_{q}).
$$
We use Lemma \ref{lem H^1(covering) = H^1} and covering $\{\mathcal U\}$ of $\mathcal B$ as above. By Theorem \ref{theor H^1(O)=0} and Theorem \ref{theor Kunneth Formula} we obtain that 
$$
H^1 (\pi^{-1}(\mathcal U_0), \widehat{\mathcal T}^h_{pq})=0.
$$ 
Further by Lemma \ref{lem H^0(flag)=0} and Theorem \ref{theor Kunneth Formula} we have  
$$
H^1(\{\pi^{-1}(\mathcal U) \}, \widehat{\mathcal T}^h_{pq})=\lbrace0\rbrace\quad \text{for} \quad p\ne q
$$ and 
$$
H^1(\{\pi^{-1}(\mathcal U) \}, \widehat{\mathcal T}^h_{pp})=H^1(\{\mathcal U \}, (\tilde{\mathcal T}_{\mathcal B})_p).
$$ 
Now our result follows from the corresponding result for super-grassmannians,  in \cite[Theorem 1]{onirig}.$\square$

\medskip

Let us recall some results of Bott \cite{Bott}, see also \cite[Section 3]{onirig}, that we will essentially use to prove the next theorem. Let $\mathbf E$ be a homogeneous bundle over $\mathcal B_0$ and $x_0\in \mathcal B_0$ be the origin of $Z_{I_1}\subset \mathcal B_0$, see (\ref{eq I_0 and I_1 of Z_I}). Denote by $\mathbf E_{x_0}$ the $P_{\mathcal B}$-module corresponding to $\mathbf E$. In other words $\mathbf E_{x_0}$ is the fiber of  $\mathbf E$ over $x_0$.   Denote also by $\mathfrak P$ and by $\mathfrak R$ the Lie algebras of $P_{\mathcal B}$ and $R_{\mathcal B}$, respectively. We have the following isomorphisms \cite[Theorem I, Corollory 2 of Theorem W$_2$]{Bott}
\begin{equation}\label{eq isom of H^1}
H^1(\mathcal B_0, \mathbf E)^{\mathfrak{gl}_{m|n}(\mathbb C)_{\bar 0}}  \simeq H^1(\mathfrak P, \mathfrak R, \mathbf E_{x_0}) \simeq H^1(\mathfrak n_{\mathcal B}, \mathbf E_{x_0})^{R_{\mathcal B}},
\end{equation}
where $\mathfrak n_{\mathcal B}$ is the nilradical of the Lie algebra of   $P_{\mathcal B}$. 
Now we are ready to prove the following theorem. 

\begin{theorem}\label{theor H^1(T_22) =0}
Assuming (\ref{eq generic type condition for flags}) and (\ref{condition flag type 1}), we have	$H^1(\mathcal M_0, \widehat{\mathcal T}_{22})=\lbrace0\rbrace$.
\end{theorem}

\smallskip

\noindent {\it Proof.}
\noindent{\bf Step 1. $H^1(\mathcal M_0,\widehat{\mathcal T}_{22})$ as invariant Lie algebra cohomology.} The following sequence of sheaves
\begin{equation}\label{eq T^v_22 to T_22 to T^h_22}
0\to \pi_*(\widehat{\mathcal T}^v_{22})\to \pi_*(\widehat{\mathcal T}_{22})\to \pi_*(\widehat{\mathcal T}^h_{22})\to 0.
\end{equation}
is exact. Indeed, from Proposition \ref{prop_geometr A,B,T} and Theorem \ref{theor H^1(T^h_pq,T^v_pq)} it follows that 
$$
H^1(\mathcal U_0,\pi_*(\widehat{\mathcal T}^v_{22}))=\lbrace0\rbrace.
$$  
Further we use the long exact sequence over any $V\subset \mathcal U_0$.  We also can choose $\mathcal U$ such that Sequence \ref{eq T^v_22 to T_22 to T^h_22} is split over $\mathcal U$. In other words we have
$$
\pi_*(\widehat{\mathcal T}_{22})\vert_{\mathcal U} = \pi_*(\widehat{\mathcal T}^v_{22})\vert_{\mathcal U}  \oplus \pi_*(\widehat{\mathcal T}^h_{22})\vert_{\mathcal U}
$$
for any $\mathcal U$.  As above we get that 
 $$
 H^1(\mathcal U,  \pi_*(\widehat{\mathcal T}^v_{22})) =0, \quad  H^1(\mathcal U,  \pi_*(\widehat{\mathcal T}^h_{22})) =0.
 $$
 Therefore we can use Lemma \ref{lem H^1(covering) = H^1} for the shaef $\widehat{\mathcal T}_{22}$. By Theorem \ref{theor H^1(T^h_pq,T^v_pq)}, Theorem \ref{theor H^1(T^h_22)=C}, the observations above and Sequence  (\ref{eq T^v_22 to T_22 to T^h_22}) we obtain
 \begin{equation}\label{eq exact seq H^1(T_22) to H^1(T^h)}
0\to  H^1(\mathcal B_0,\pi_*(\widehat{\mathcal T}_{22}))\to H^1(\mathcal B_0,\pi_*(\widehat{\mathcal T}^h_{22}))= \mathbb C. 
 \end{equation}

Denote by $\mathbf T$, $\mathbf T^v$ and by  $\mathbf T^h$ the vector bundles corresponding to the locally free sheaves $\pi_*(\widehat{\mathcal T}_{22})$, $\pi_*(\widehat{\mathcal T}^v_{22})$ and $\pi_*(\widehat{\mathcal T}^h_{22})$ on $\mathcal B_0$, respectively. Denote also by $\mathbf T_{\mathcal B}$ the vector bundle corresponding to the locally free sheaf $(\tilde{\mathcal T}_{\mathcal B})_2$ on $\mathcal B_0$. Further using (\ref{eq exact seq H^1(T_22) to H^1(T^h)}) and (\ref{eq isom of H^1}) we get 
$$
H^1(\mathcal B_0,\pi_*(\widehat{\mathcal T}_{22}))= H^1(\mathcal B_0,\pi_*(\widehat{\mathcal T}_{22}))^{\mathfrak{gl}_{m|n}(\mathbb C)_{\bar 0}}= H^1(\mathfrak n_{\mathcal B},\mathbf T_{x_0})^{R_{\mathcal B}}.
$$

\noindent{\bf Step 2. Computation of  $H^1(\mathfrak n_{\mathcal B},\mathbf T_{x_0})^{R_{\mathcal B}}$.} 
 A direct computation shows that there is the following $R_{\mathcal B}$-invariant decomposition 
$$
\mathbf T_{x_0}= \mathbf T^v_{x_0}\oplus \mathbf T^h_{x_0}= \mathbf T^v_{x_0}\oplus (\mathbf T_{\mathcal B})_{x_0}.
$$ 
 Consider the chart corresponding to the coordinate matrices $Z_{I_i}$ given by (\ref{eq I_0 and I_1 of Z_I}). Let us write $Z_{I_s}$ explicitly
\begin{equation}
\label{Z_I_s explicit}
\begin{split}
Z_{I_s}=\left(
\begin{array}{cc}
X_s & \Xi_s\\
E & 0\\
\H_s & Y_s \\
0 & E
\end{array} \right), \ \ s=1,\dots,r.,
\end{split}
\end{equation}
where $X_s= (x^s_{ij})$, $Y_s= (y^s_{ij})$, $\Xi_s= (\xi^s_{ij})$ and $\H_s= (\eta^s_{ij})$. In \cite[Proof of Theorem 1]{onirig} a basis of vector space $C^1(\mathfrak n_{\mathcal B},(\mathbf T_{\mathcal B})_{x_0})^{R_{\mathcal B}}$ of cochains was calculated.
 This basis contains only the following two elements 
$$
c(e_{\alpha\beta})= a\sum_{ij} \eta^1_{i\alpha}\xi^1_{\beta j}\frac{\partial}{\partial y^1_{ij}};\quad
c(f_{\alpha\beta})= b \sum_{ij} \xi^1_{i\alpha} \eta^1_{\beta j}\frac{\partial}{\partial x^1_{ij}}, \quad a,b\in \mathbb C.
$$
Here $e_{\alpha\beta}\in \mathfrak n^1_{\mathcal B}\simeq \Mat_{k_1, m-k_1}(\mathbb C)$, $e_{\alpha\beta}\in \mathfrak n^2_{\mathcal B}\simeq \Mat_{l_1, n-l_1}(\mathbb C)$ and $\mathfrak n_{\mathcal B}= \mathfrak n^1_{\mathcal B}\oplus \mathfrak n^2_{\mathcal B}$, where $\mathfrak n^1_{\mathcal B}$ and $\mathfrak n^2_{\mathcal B}$ are irreducible $R_{\mathcal B}$-modules with highest weights $-\mu_{m-k_1}+\mu_{m-k_1+1}$ and $ - \lambda_{n-l_1}+\lambda_{n-l_1+1}$, respectively. From the proof of Theorem \ref{theor A_p= and A_pq =}, Step 2, and from the proof of Theorem \ref{theor H^1(B^h_pq)=}, Step 2,  it follows that the representation of $R_{\mathcal B}$ in  $\mathbf T^v_{x_0}$ does not have weights of index $1$. In particular it follows that this representation does not have weights $-\mu_{m-k_1}+\mu_{m-k_1+1}$ and $ - \lambda_{n-l_1}+\lambda_{n-l_1+1}$. Therefore any cochain $d\in C^1(\mathfrak n_{\mathcal B},\mathbf T_{x_0})^{R_{\mathcal B}}$ has values in $(\mathbf T_{\mathcal B})_{x_0}$. Therefore we have $C^1(\mathfrak n_{\mathcal B},\mathbf T_{x_0})^{R_{\mathcal B}} \simeq  C^1(\mathfrak n_{\mathcal B},(\mathbf T_{\mathcal B})_{x_0})^{R_{\mathcal B}}$.

Let us show that from $\delta(c)=0$ for $c\in C^1(\mathfrak n_{\mathcal B},\mathbf T_{x_0})^{R_{\mathcal B}}$ it follows that $c=0$. In other words let us show that the vector space of cocycles $Z^1(\mathfrak n_{\mathcal B},\mathbf T_{x_0})^{R_{\mathcal B}}$ is trivial. By definition $c\in Z^1(\mathfrak n_{\mathcal B},\mathbf T_{x_0})^{R_{\mathcal B}}$ if and only if
\begin{equation}\label{eq d(c)=0}
(\delta c)(x,y)=xc(y)-yc(x)=0 \,\,\,\text{for any}\,\,\, x,y\in \mathfrak n_{\mathcal B}.
\end{equation}
Let us calsulate the fundamental vector fields on $\widehat{\mathcal M}$ corresponding to the matrices $U_1\in \mathfrak n^1_{\mathcal B}$ and $U_2\in \mathfrak n^2_{\mathcal B}$ defined by 
$$
\left(
\begin{array}{cccccc}
E&0&0&0&0&0\\
U_1&E&0&0&0&0\\
0&0&E&0&0&0\\
0&0&0&E&0&0\\
0&0&0&U_2&E&0\\
0&0&0&0&0&E
\end{array}
\right
)
\left(
\begin{array}{cccc}
X^{11}&X^{12}&\Xi^{11}&\Xi^{12}\\
E&0&0&0\\
0&E&0&0\\
\H^{11}&\H^{12}&Y^{11}&Y^{12}\\
0&0&E&0\\
0&0&0&E
\end{array}
\right
)
\left(
\begin{array}{cc}
X^2&\Xi^2\\
E&0\\
\H^2&Y^2\\
0&E
\end{array}
\right)\cdots .
$$
Here $U_1\in \Mat_{k_1-k_2,m-k_1}(\mathbb C)$ and $U_2\in \Mat_{l_1-l_2,n-l_1}(\mathbb C)$. The first matrix is divided in blocks of size $(m-k_1,k_1-k_2,k_2, n-l_1,l_1-l_2,l_2)$, other matrices are divided in blocks accordingly. We will denote by $\tilde U_i^*$ the fundamental vector fields on $\tilde{\mathcal M}$ corresponding to $U_i$. Consider (\ref{eq d(c)=0}). As we have seen above $c(x)\in (\mathbf T_{\mathcal B})_{x_0}$. Further any fundamental vector field is projectible. Therefore we can decompose (\ref{eq d(c)=0}) into the horizontal and the vertical parts. Our idea is to show that the vertical part of (\ref{eq d(c)=0}) does not hold true.

After a direct computation we get $\tilde U_i^*= (\tilde U_i^*)^h + (\tilde U_i^*)^v $, where 
\begin{align*}
(\tilde U_1^*)^v=  &
\sum_{\alpha\beta}( U_1X^{11}X^2+ U_1X^{12})_{\alpha\beta}\frac{\partial}{\partial x^{2}_{\alpha\beta}}+\\
&\sum_{\alpha\beta}( U_1X^{11}\Xi^2+ U_1\Xi^{11}Y^2+ U_1\Xi^{12})_{\alpha\beta}\frac{\partial}{\partial \xi^{2}_{\alpha\beta}};\\
(\tilde U_2^*)^v= &\sum_{\alpha\beta}( U_2Y^{11}Y^2+ U_2Y^{12})_{\alpha\beta}\frac{\partial}{\partial y^{2}_{\alpha\beta}} \\
&+
\sum_{\alpha\beta}( U_2Y^{11}\H^2+ U_2\H^{11}X^2+ U_2\H^{12})_{\alpha\beta}\frac{\partial}{\partial \eta^{2}_{\alpha\beta}}.
\end{align*}
Denote by $\widehat U_i^*$ the fundamental vector field on $\widehat{\mathcal M}$ corresponding to $U_i$. 
We have $\widehat U_i^*= (\widehat U_i^*)^h + (\widehat U_i^*)^v $, where
\begin{align*}
(\widehat U_1^*)^v=  &
\sum_{\alpha\beta}( U_1X^{11}X^2+ U_1X^{12})_{\alpha\beta}\frac{\partial}{\partial x^{2}_{\alpha\beta}}+ \sum_{\alpha\beta}( U_1X^{11}\Xi^2)_{\alpha\beta}\frac{\partial}{\partial \xi^{2}_{\alpha\beta}};\\
(\widehat U_2^*)^v= &\sum_{\alpha\beta}( U_2Y^{11}Y^2+ U_2Y^{12})_{\alpha\beta}\frac{\partial}{\partial y^{2}_{\alpha\beta}} +
\sum_{\alpha\beta}(U_2Y^{11}\H^2)_{\alpha\beta}\frac{\partial}{\partial \eta^{2}_{\alpha\beta}}.
\end{align*}

Therefore
\begin{align*}
(\widehat e^*_{ij})^v= \sum_{\alpha\beta}x^{11}_{j\alpha}x^2_{\alpha\beta}\frac{\partial}{\partial x^{2}_{i\beta}}+
\sum_{\alpha}x^{12}_{j\alpha}\frac{\partial}{\partial x^{2}_{i\alpha}}+
\sum_{\alpha\beta}x^{11}_{j\alpha}\xi^2_{\alpha\beta}\frac{\partial}{\partial \xi^{2}_{i\beta}}
;\\
(\widehat f^*_{kl})^v= \sum_{\alpha\beta}y^{11}_{l\alpha}y^2_{\alpha\beta}\frac{\partial}{\partial y^{2}_{k\beta}}+
\sum_{\alpha}y^{12}_{l\alpha}\frac{\partial}{\partial y^{2}_{k\alpha}}+
\sum_{\alpha\beta}y^{11}_{l\alpha}\eta^2_{\alpha\beta}\frac{\partial}{\partial \eta^{2}_{k\beta}}.
\end{align*}

Let us compute the vertical part of $(\delta c)$. We have

\begin{align*}
(\widehat e^*_{ij})^vc(f_{kl})=& -b (\sum_{\alpha\beta} \xi^{11}_{j k}\eta^{11}_{l \alpha}x^2_{\alpha\beta}\frac{\partial}{\partial x^{2}_{i\beta}}+
\sum_{\alpha} \xi^{11}_{j k}\eta^{12}_{l \alpha}\frac{\partial}{\partial x^{2}_{i\alpha}}+
\sum_{\alpha\beta} \xi^{11}_{j k}\eta^{11}_{l \alpha}\xi^2_{\alpha\beta}\frac{\partial}{\partial \xi^{2}_{i\beta}}).\\
(\widehat f^*_{kl})^v c(e_{ij})= &-a( \sum_{\alpha\beta} \eta^{11}_{l i}\xi^{11}_{j \alpha}x^2_{\alpha\beta}\frac{\partial}{\partial y^{2}_{k\beta}}+
\sum_{\alpha} \eta^{11}_{l i}\xi^{12}_{j \alpha}\frac{\partial}{\partial y^{2}_{k\alpha}}+
\sum_{\alpha\beta} \eta^{11}_{l i}\xi^{11}_{j \alpha}\eta^2_{\alpha\beta}\frac{\partial}{\partial \eta^{2}_{k\beta}}).
\end{align*}
From (\ref{eq d(c)=0}) we get $a=b=0$, therefore $Z^1(\mathfrak n_{\mathcal B},\mathbf T_{x_0})^{R_{\mathcal B}}= \lbrace0\rbrace$, hence 
$$
H^1(\mathcal M_0,\widehat{\mathcal T}_{22})=  H^1(\mathfrak n_{\mathcal B},\mathbf T_{x_0})^{R_{\mathcal B}}=\lbrace0\rbrace.
$$ 
The proof is complete$\square$

\subsection{1st cohomology group with values in the sheaf $\widehat{\mathcal T}_{20}$}

We start with the following lemma.

\begin{lemma}\label{lem H^1(T_20)=H^1(T_20)^G}
	Assuming (\ref{eq generic type condition for flags}) and (\ref{condition flag type 1}), we have 
	$$
	H^1(\mathcal M_0, \mathcal C^v_{20})=H^1(\mathcal M_0,
	\mathcal C^v_{20})^{\mathfrak{gl}_{m|n}(\mathbb C)_{\bar 0}}\quad \text{and} \quad H^1(\mathcal M_0, \widehat{\mathcal T}^v_{20})=H^1(\mathcal M_0, \widehat{\mathcal T}^v_{20})^{\mathfrak{gl}_{m|n}(\mathbb C)_{\bar 0}}.
	$$
\end{lemma}

\noindent {\it Proof.} The second statement follows from the first one, from Theorem \ref{theor A_p= and A_pq =}, Theorem \ref{theor H^1(B^h_pq)=} and from the following exact sequence
$$
0\to H^1(\mathcal M_0, \mathcal A^v_{20})= \lbrace0\rbrace\to H^1(\mathcal M_0, \widehat{\mathcal T}^v_{20})\to H^1(\mathcal M_0, \mathcal C^v_{20}).
$$

Let us prove the first statement. Consider the representation $\tau^v\otimes \bigwedge^2\psi$ of $P$, see Lemma \ref{lemma_reps of O_11 nd O_10} and Lemma \ref{lemma_reps of T,T^v and T/T^v},   corresponding to the homogeneous locally free  sheaf $\mathcal C^v_{20}= \Theta^v\otimes \bigwedge^2 \mathcal O_{10}$. A weight of this representation has the form $\Lambda= \Lambda_0 + \Lambda_1$, where
$$
\Lambda_0=
\left[
\begin{array}{l}
1.\,\, \mu_i^p-\mu_j^q,\,\, 1<p<q;\\
2.\,\, 0;
\end{array}
\right.
+
\left[
\begin{array}{l}
A.\,\, -\mu_{i_1}^s-\mu_{i_2}^t,\,\, 1<s,t<r+1;\\
B.\,\, \mu_{j_1}^s+\mu_{j_2}^t,\,\, s,t>2;\\
C.\,\, -\mu_{i_1}^s+\mu_{j_1}^t,\,\, t>2, 1<s<r+1;\\
\end{array}
\right.
$$
$$
\Lambda_1=
\left[
\begin{array}{l}
1.\,\, 0;\\
2.\,\, \lambda_i^p-\lambda_j^q,\,\, 1<p<q;\\
\end{array}
\right.
+
\left[
\begin{array}{l}
A.\,\, \lambda_{j_1}^s+\lambda_{j_2}^t,\,\, s,t>2;\\
B.\,\, -\lambda_{i_1}^s-\lambda_{i_2}^t,\,\, 1<s,t<r+1;\\
C.\,\, -\lambda_{i_2}^s+\lambda_{j_2}^t,\,\, t>2, 1<s<r+1.\\
\end{array}
\right.
$$
Assume that $\Lambda_0$ is dominant and $\Lambda_1$ has index $1$. In case $(1.A)$ the weight $\Lambda_1$ is singular. In case $(1.B)$ the weight $\Lambda_0$ cannot be dominant. In case $(1.C)$ under our assumptions we have only one possibility $\Lambda_0=0$ and $\Lambda_1= -\lambda^{s}_{i_1}+\lambda^{s+1}_{i_1+1}$. It is easy to see that the reflection $\mathrm{s}$ corresponding to the root $\lambda^{s}_{i_1} -  \lambda^{s+1}_{i_1+1}$ maps $\Lambda + \zeta$, where $\zeta$ is as in  (\ref{eq half-sum of positive roots}), to the weight $\zeta$. Since $\mathrm{s}(\Lambda + \zeta) - \zeta =0$. In other words the weight $\Lambda + \zeta$ corresponds to $1$-dimensional trivial $\mathfrak{gl}_{m|n}(\mathbb C)_{\bar 0}$-module.

In cases $(2.A)$ and $(2.B)$ the weight  $\Lambda_0$ is not dominant. In case $(2.C)$ if $\Lambda_0$ is dominant we have only one possibility $\Lambda_0=0$. For $\Lambda_1$ we also have only one possibility $\Lambda_1= -\lambda_{i}+\lambda_{i+1}$ for some $i$. This weight also  corresponds to $1$-dimensional trivial $\mathfrak{gl}_{m|n}(\mathbb C)_{\bar 0}$-module.$\square$

\medskip

The following theorem was proved in \cite{onirig}.

\begin{theorem}\cite[Theorem 1]{onirig} Assume that (\ref{eq generic type condition for flags}) and (\ref{condition flag type 1}) hold true and $r=1$. Then 
	$$
	 H^1(\mathcal M_0, \tilde{\mathcal T}_{2}) = \mathbb C.
	$$
	
\end{theorem}

\medskip

Now we are ready to prove the following theorem.

\begin{theorem}\label{theor H^1(T^v_20)=C} 
	Assume that (\ref{eq generic type condition for flags}) and (\ref{condition flag type 1}) hold true. If  $H^1(\mathcal S_0, (\tilde{\mathcal T}_{\mathcal S})_{2})= \mathbb C$, then 
 $$
 H^1(\mathcal M_0, \widehat{\mathcal T}^v_{20}) = \mathbb C.
 $$
\end{theorem}

\smallskip

\noindent {\it Proof.} Assume that $H^1(\mathcal S_0, (\tilde{\mathcal T}_{\mathcal S})_{2})= \mathbb C$.  By Lemma \ref{lem H^1(T_20)=H^1(T_20)^G} we need to calculate $H^1(\mathcal M_0, \mathcal T^v_{20})^{\mathfrak{gl}_{m|n}(\mathbb C)_{\bar 0}}$. Denote by $P_{\mathcal S}$ a copy of the group $P$ for the supermanifold $\mathcal S$ and by $R_{\mathcal S}$ the reductive part of $P_{\mathcal S}$. Denote also by $\mathbf T^v_{20}$ and $(\mathbf T_{\mathcal S})_2$ the vector bundles corresponding to the locally free sheaves $\widehat{\mathcal T}^v_{20}$ and $(\tilde{\mathcal T}_{\mathcal S})_{2}$.  Let us prove the equality 
$$
H^1(\mathfrak n, (\mathbf T^v_{20})_{x_0})^R= H^1(\mathfrak n_{\mathcal S}, ((\mathbf T_{\mathcal S})_{2})_{\tilde x_0})^{R_{\mathcal S}},
$$ 
where $x_0=P$,  $\tilde x_0= P_{\mathcal S}$ and $\mathfrak n_{\mathcal S}$ is the nilradical of the Lie algebra of  $P_{\mathcal S}$. As usual we denote by $\mathbf E_x$ the fiber of a vector bundle $\mathbf E$.  Note that by Proposition \ref{prop_geometr A,B,T} we have  
$$
(\mathbf T^v_{20})_{x_0}= ((\mathbf T_{\mathcal S})_{2})_{\tilde x_0}.
$$ 
The representation of $P_{\mathcal S}$ in $((\mathbf T_{\mathcal S})_{2})_{\tilde x_0}$ is equal to the restriction of the representation of $P$ in  $(\mathbf T^v_{20})_{x_0}$ to the subgroup $P_{\mathcal S}\subset P$. The $R$-module $\mathfrak n$ is a completely reducible $R$-module with highest weights 
$$
-\mu^p_{max}+\mu^q_{min},\quad  -\lambda^p_{max}+\lambda^q_{min}, \quad  p<q, \quad p,q\in \lbrace 1, \ldots, r+1\rbrace.
$$
The $R_{\mathcal S}$-module  $\mathfrak n_{\mathcal S}$ is a completely reducible $R_{\mathcal S}$-module with highest weights 
$$
-\mu^p_{max}+\mu^q_{min},\quad -\lambda^p_{max}+\lambda^q_{min}, \quad p<q, \quad  p,q\in \lbrace 2, \ldots, r+1\rbrace.
$$ 
We denoted by $\mu^p_{max}$ the element $\mu^p_{i}$ with maximal possible index $i$ from the block $p$ and by $\mu^q_{min}$ the element $\mu^p_{j}$ with minimal possible index $j$ from the block $q$. The same agreement we use for $\lambda$'s.

Further we note that the weights of $P$ in a fiber of $(\mathbf T^v_{20})_{x_0}$ is equal to the weights of $P$ in the fiber of the sheaf $\mathcal A^v_{20}$ plus  the weights of $P$ in the fiber of the sheaf $\mathcal C^v_{20}$. Above we calculated the representations corresponding to the homogeneous sheaves $\mathcal A^v_{20}$ and $\mathcal C^v_{20}$. Using Lemma \ref{lemma_reps of O_11 nd O_10} and Lemma \ref{lemma_reps of T,T^v and T/T^v}, we get that they have the following form 
$$
\phi^*\otimes \bigwedge^3\phi\quad  \text{and} \quad \tau^v\otimes \bigwedge^2\phi,
$$ 
respectively. We see that the corresponding weights do not contain elements $-\mu^1_{max}$ and $-\lambda^1_{max}$. Moreover we have the following decomposition of $R$-modules $\mathfrak n= \mathfrak n_{\mathcal B}\oplus \mathfrak n_{\mathcal S}$, where $\mathfrak n_{\mathcal B}$ is a completely reducible $R$-module with highest weights $-\mu^1_{max}+\mu^q_{min}$ and $-\lambda^1_{max}+\lambda^q_{min}$, where $q>1$. (Compare also which the proof of Theorem \ref{theor H^1(T_22) =0}.) Now we see that  we have the equality of vector spaces of cochains 
$$
C^1(\mathfrak n, (\mathbf T^v_{20})_{x_0})^R= C^1(\mathfrak n_{\mathcal S}, ((\mathbf T_{\mathcal S})_{2})_{\tilde x_0})^{R_{\mathcal S}}.
$$

Let us prove that $\mathfrak n_{\mathcal B}$ acts trivially in  $(\mathbf T^v_{20})_{x_0}$. Any fundamental vector field on $\widehat{\mathcal M}$ has the form  $v=v^h+v^v$, where $v^h\in \widehat{\mathcal T}^h_{00}$ and  $v^v\in \widehat{\mathcal T}^v_{00}$. The part $v^h$ acts on  $(\mathbf T^v_{20})_{x_0}$ trivially since sections of the sheaf $\widehat{\mathcal T}^v_{20}$ do not depend on coordinates $Z_{I_1}$. Further the part $v^v$ also acts trivially on $(\mathbf T^v_{20})_{x_0}$. Indeed, $\mathfrak n_{\mathcal B}$ contains all matrices of the following form
$$
\left(
\begin{array}{cc}
0& 0\\
C_1& 0
\end{array}
\right) 
\times
\left(
\begin{array}{cc}
0& 0\\
C_2& 0
\end{array}
\right), 
$$
where $C_1\in \Mat_{k_1\times {m-k_1}}(\mathbb C)$ and $C_2\in \Mat_{l_1\times {n-l_1}}(\mathbb C)$. In coordinates $(Z_{I})$, see (\ref{eq I_0 and I_1 of Z_I}) for definition of $I$, after a short calculation we get 
\begin{equation}\label{eq action on S}
\left(
\begin{array}{cc}
E+C_1X^1& C_1\Xi^1\\
C_2\H^1& E+C_2Y^1
\end{array}
\right)\,\,Z_{I_2}, Z_{i_3}, \cdots .
\end{equation}
To get the action in  $(\mathbf T^v_{20})_{x_0}$  we need to put $\Xi^1=\H^1=0$ and $X^1= Y^1=0$ in  (\ref{eq action on S}). Now  we see that this action is trivial. Therefore the action of $\mathfrak n$ and of  $\mathfrak n_{\mathcal S}$ in the fibers $(\mathbf T^v_{20})_{x_0}$ and $((\mathbf T_{\mathcal S})_{2})_{\tilde x_0}$ coincide.  Therefore, 
$$
H^1(\mathfrak n, (\mathbf T^v_{20})_{x_0})^R= H^1(\mathfrak n_{\mathcal S}, ((\mathbf T_{\mathcal S})_{2})_{\tilde x_0})^{R_{\mathcal S}},
$$ 
The result follows.$\square$

\medskip

We will use the notation $V \leq\mathbb C$ for a vector space $V$ over $\mathbb C$. This meant that $V$ is either equal to $\mathbb C$ or $\{0\}$. 
From Theorem \ref{theor H^1(T^v_20)=C} we get. 

\begin{corollary}
	Assume (\ref{eq generic type condition for flags}) and (\ref{condition flag type 1}). If  $H^1(\mathcal S_0, (\tilde{\mathcal T}_{\mathcal S})_{2})= \mathbb C$, then  
	$
	H^1(\mathcal M_0, \widehat{\mathcal T}_{20})\leq\mathbb C.
	$
\end{corollary}

\smallskip

\noindent {\it Proof.} Consider the following exact sequence 
$$
0\to \mathcal T^v_{20}\to \mathcal T_{20}\to \mathcal T^h_{20}\to 0.
$$
Therefore we get
$$
\mathbb C\to H^1(\mathcal M_0, \widehat{\mathcal T}_{20})\to 0. \quad \square
$$

We are ready to prove one of our main results. 

\begin{theorem}\label{theor H^1(T_2)=C}
	Assuming (\ref{eq generic type condition for flags}) and (\ref{condition flag type 1}), we have 
	$$
	H^1(\mathcal M_0, \tilde{\mathcal T}_{2})=\mathbb C
	$$
\end{theorem}

\noindent {\it Proof.} Under Theorem's assumption the flag supermanifold $\mathcal M$ is not split. This fact can be deduced for instance from results \cite{ViGL in JA}, where the Lie algebras of holomorphic vector fields were computed. The idea is the following:  a supermanifold is split if and only if it possesses a so called grading vector field.  (Details about grading vector fields can be found in \cite{ViSplitting problem}.)

Further in \cite{Green} it was shown that the classes of non-split supermanifolds are in one-to-one correspondence with orbits of a certain group acting on $H^1(\mathcal M_0, \tilde{\mathcal T}_{2})$, and the point $0$ corresponds to the split supermanifold.  Therefore from the fact that a flag supermanifold is not split it follows that $H^1(\mathcal M_0, \tilde{\mathcal T}_{2})$ is not zero. Now consider the following exact sequence of sheaves
$$
0\to \tilde{\mathcal T}_{2(q+1)}\to \tilde{\mathcal T}_{2(q)}\to \widehat{\mathcal T}_{2q}\to 0.
$$
For $q=3$ we have an isomorphism 
$\tilde{\mathcal T}_{2(3)}\simeq \widehat{\mathcal T}_{23}$. Therefore, 
$$
H^1(\mathcal M_0, \tilde{\mathcal T}_{2(3)})= H^1(\mathcal M_0, \widehat{\mathcal T}_{23})= \lbrace0\rbrace.
$$ 
Further
$$
0\to H^1(\mathcal M_0, \tilde{\mathcal T}_{2(2)})\to H^1(\mathcal M_0, \widehat{\mathcal T}_{22})= \lbrace0\rbrace,
$$
hence, $H^1(\mathcal M_0, \tilde{\mathcal T}_{2(2)})= \lbrace0\rbrace$. Similarly, $H^1(\mathcal M_0, \tilde{\mathcal T}_{2(1)})= \lbrace0\rbrace$. From the exact sequence 
$$
0\to H^1(\mathcal M_0, \tilde{\mathcal T}_{2(0)})\to H^1(\mathcal M_0, \widehat{\mathcal T}_{20})\leq \mathbb C
$$
it follows that $H^1(\mathcal M_0, \mathcal T_{2(0)})\leq \mathbb C$. Finally we obtain
$$
\leq\mathbb C\to H^1(\mathcal M_0, \tilde{\mathcal T}_{2(-1)})\to H^1(\mathcal M_0, \widehat{\mathcal T}_{20})= \lbrace0\rbrace.
$$
Therefore, $H^1(\mathcal M_0, \tilde{\mathcal T}_{2})=H^1(\mathcal M_0, \tilde{\mathcal T}_{2(-1)})\leq \mathbb C$ $\square$

\smallskip
\begin{corollary}
	Under the conditions of Theorem \ref{theor H^1(T_2)=C}, the flag supermanifold $\mathcal M$ is the unique non-split supermanifold up to isomorphism which corresponds to the split supermanifold $\tilde{\mathcal M}$. 
\end{corollary}

Now we are ready to prove another main result. 

\begin{theorem}\label{theor main H^1(T)=0} Assuming (\ref{eq generic type condition for flags}) and (\ref{condition flag type 1}), we have 
	$H^1(\mathcal M_0, \mathcal T)= \lbrace0\rbrace$.
\end{theorem}

\smallskip

\noindent {\it Proof.} The following sequence
$$
0\to \mathcal T_{(p+1)}\to \mathcal T_{(p)}\to \mathcal T_{p}\to 0
$$
is exact. Hence, $H^1(\mathcal M_0,\mathcal T_{(p)})= \lbrace0\rbrace$ for $p\leq 3$ and $H^1(\mathcal M_0,\mathcal T_{(2)})\leq \mathbb C$. From $H^0(\mathcal M_0,\tilde{\mathcal T}_{1})= \lbrace0\rbrace$ and from $H^1(\mathcal M_0,\tilde{\mathcal T}_{1})= \lbrace0\rbrace$ we get $H^1(\mathcal M_0,\mathcal T_{(1)})\simeq H^1(\mathcal M_0,\mathcal T_{(2)})$. Further, for $p=0$ we have
\begin{align*}
0\to &H^0(\mathcal M_0,{\mathcal T}_{(1)}) \to H^0(\mathcal M_0,{\mathcal T}_{(0)}) \to H^0(\mathcal M_0,\tilde{\mathcal T}_{0}) \to\\
& H^1(\mathcal M_0,{\mathcal T}_{(1)}) \to H^1(\mathcal M_0,{\mathcal T}_{(0)}) \to H^1(\mathcal M_0,\tilde{\mathcal T}_{0})
\end{align*}
Under our assumption $H^0(\mathcal M_0,{\mathcal T}_{(0)}) \simeq \mathfrak{pgl}_{m|n}(\mathbb C)$. (This result was obtained in \cite{ViGL in JA}). Further by Theorem \ref{theor vect fields on flag}, $H^0(\mathcal M_0,\tilde{\mathcal T}_{0}) \simeq \mathfrak{gl}_{m|n}(\mathbb C)$. The image 
$H^0(\mathcal M_0,{\mathcal T}_{(0)}) \to H^0(\mathcal M_0,\tilde{\mathcal T}_{0})$ has codimendsion $1$ and contains the grading operator.  We also proved that $H^1(\mathcal M_0,\tilde{\mathcal T}_{0}) = \{0\}$. From these facts we conclude that $H^1(\mathcal M_0,{\mathcal T}_{(0)})=\{0\}$. Therefore,
$$
0 \to H^1(\mathcal M_0,{\mathcal T}_{(-1)}) \to H^1(\mathcal M_0,\tilde{\mathcal T}_{-1}) = 0.
$$
The result follows.$\Box$

\bigskip

\noindent
E.~V.: Departamento de Matem{\'a}tica, Instituto de Ci{\^e}ncias Exatas,
Universidade Federal de Minas Gerais,
Av. Ant{\^o}nio Carlos, 6627, CEP: 31270-901, Belo Horizonte,
Minas Gerais, BRAZIL, and Laboratory of Theoretical and Mathematical Physics, Tomsk State University, 
Tomsk 634050, RUSSIA, 

\noindent email: {\tt VishnyakovaE\symbol{64}googlemail.com}

\end{document}